\newif\ifieee  
\ieeetrue
\newif\ifieeefinal
\ieeefinaltrue
\newif\ifedit
\editfalse
\newif\ifwarn   
\warntrue
\newif\ifwarnshort  
\warnshortfalse
\newif\iftodo   
\todotrue
\newif\iftodoshort  
\todoshortfalse
\newif\ifmore   
\moretrue
\newif\ifmoreshort  
\moreshorttrue
\ifieee
  \ifieeefinal
    \documentclass[letterpaper,twoside,twocolumn]{IEEEtran}
  \else
    \documentclass[letterpaper,onecolumn,draftcls]{IEEEtran}
  \fi 
\else
  \documentclass[11pt,fleqn]{article}
\fi 
\ifieee
\else
\fi
\ifedit
  \usepackage[german,english]{babel}   
\fi
\usepackage{amsmath}
\usepackage{amssymb}
\ifieee
  \usepackage{theorem}
\else
  \usepackage{amsthm}
\fi
\usepackage{graphicx}
\usepackage{epsfig}
\usepackage{xspace}
\ifedit
  \usepackage{mycolor}
  \usepackage{manfnt}  
  \usepackage{marvosym}  
  \DeclareMathSymbol{\myRightarrow}{3}{symbols}{41} 
\fi
\newcommand{\myline}{\hfill\\} 
\newcommand{\emptypage}%
{
  \newpage
  \vspace*{10cm}
  \pagebreak
}
\ifieee
  \newcommand{\textin}{\hangindent0.5cm\hangafter1\myline}
  \newcommand{\textout}{\hangindent0cm\hangafter1\myline}
\else
  \newcommand{\textin}{\hangindent0.75cm\hangafter1\myline}
  \newcommand{\textout}{\hangindent0cm\hangafter1\myline}
\fi
\newcommand{\more}[2]     
{\ifmore
  \ifmoreshort
    \reversemarginpar\marginpar{{\sf\footnotesize More?}\,
      \begin{minipage}[c]{0.05\linewidth}\tiny #1 \end{minipage}}
  \else
    \par
    \reversemarginpar\marginpar{{\sf\footnotesize More?}\,
      \begin{minipage}[c]{0.05\linewidth}\tiny #1 \end{minipage}}
    \fbox{\begin{minipage}[t]{0.925\linewidth}\selectlanguage{german}
       #2
    \end{minipage}
    }
    \par
  \fi 
\fi}
\ifieee
  \theoremstyle{plain}
  \theoremheaderfont{\bf}
  {\theorembodyfont{\it}
  \newtheorem{theorem}{Theorem}[section] 
  \newtheorem{proposition}[theorem]{Proposition}
  \newtheorem{lemma}[theorem]{Lemma}
  \newtheorem{corollary}[theorem]{Corollary}
  \newtheorem{conclusion}[theorem]{Conclusion}
  }
  {\theorembodyfont{\rm}
  \newtheorem{definition}[theorem]{Definition}
  }
  {\theorembodyfont{\rm}
  \newtheorem{remark}[theorem]{Remark}
  \newtheorem{note}[theorem]{Note}
  \newtheorem{cit}[theorem]{Citation}
  }
\else
  \theoremstyle{plain}
  \newtheorem{theorem}{Theorem}[section]  
  \newtheorem*{theorem*}{Theorem}
  \newtheorem{proposition}[theorem]{Proposition}
  \newtheorem*{proposition*}{Proposition}
  \newtheorem{lemma}[theorem]{Lemma}
  \newtheorem*{lemma*}{Lemma}
  \newtheorem{corollary}[theorem]{Corollary}
  \newtheorem*{corollary*}{Corollary}
  
  \newtheorem*{conclusion*}{Conclusion}
  \theoremstyle{definition}
  
  \newtheorem*{definition*}{Definition}
  \theoremstyle{remark}
  \newtheorem{remark}[theorem]{Remark}
  \newtheorem*{remark*}{Remark}
  
  \newtheorem*{note*}{Note}
  
  \newtheorem*{cit*}{Citation}

\fi 
\ifieee\else
  \renewcommand{\theequation}{\arabic{equation}}
\fi
\newcommand{\mathcomma}{\:,\:\:}
\newcommand{\mapset}[2]{#1\longrightarrow #2}
\newcommand{\mapto}[2]{#1\longmapsto #2}
\usepackage{dsfont}
\newcommand{\CC}{\ensuremath{\mathds{C}}}   
\newcommand{\RR}{\ensuremath{\mathds{R}}}   
   
\newcommand{\ZZ}{\ensuremath{\mathds{Z}}}   
\newcommand{\NN}{\ensuremath{\mathds{N}}}

\newcommand{\mdef}{:=}
\newcommand{\mfed}{=:}
\newcommand{\oneover}[1]{\frac{1}{#1}}
\newcommand{\toneover}[1]{\tfrac{1}{#1}}
\newcommand{\comp}{\raisebox{0.4ex}{$\scriptscriptstyle\circ\ $}}
\renewcommand{\d}{\partial} 
\newcommand{\norm}[1]{\lVert#1\rVert}

\newcommand{\abs}[1]{\lvert#1\rvert}
\newcommand{\bigabs}[1]{\left\lvert#1\right\rvert}
\newcommand{\generate}[1]{\left\langle #1 \right\rangle}

\newcommand{\scalprod}[1]{<\!\!#1\!\!>}

\newcommand{\lnorm}[2]{\lVert#1\rVert_{\scriptscriptstyle L^{#2}}}

\newcommand{\fnorm}[1]{\lVert#1\rVert_{\scriptscriptstyle\text{F}}}

\newcommand{\Unity}{\ensuremath{\mathbf{1}}\xspace}
\newcommand{\Zero}{\ensuremath{\mathbf{0}}\xspace}

\renewcommand{\o}{\ensuremath{o}}
\newcommand{\obar}[1]{\smash[t]{\overset{\rule[-0.5pt]{0.75ex}{0.25pt}}{#1}}}
\newcommand{\ubar}[1]{\smash[b]{\underset{\rule[5pt]{0.75ex}{0.25pt}}{#1}}}

\newcommand{\textmatrix}[1]{
   \left(\begin{smallmatrix}#1\end{smallmatrix}\right)}
\newcommand{\eqmatrix}[1]{\begin{pmatrix}#1\end{pmatrix}}
\newcommand{\smallperp}{\scriptscriptstyle\perp}
\newcommand{\smallparallel}{\raisebox{0.25ex}{$\scriptscriptstyle\parallel$}}

\newcommand{\followshort}{\ensuremath{\;\Longrightarrow\;}}

\newcommand{\tends}{\ensuremath{\longrightarrow}}
\newcommand{\tendshort}{\ensuremath{\rightarrow}}
\DeclareMathOperator{\tr}{tr}

\DeclareMathOperator{\vol}{vol}

\DeclareMathOperator{\re}{\mathfrak{Re}}

\DeclareMathOperator{\Gl}{GL}

\DeclareMathOperator{\diag}{diag}

\DeclareMathOperator{\dimc}{\dim_{\CC}}
\DeclareMathOperator{\dimr}{\dim_{\RR}}

\DeclareMathOperator{\ric}{Ric}
\DeclareMathOperator{\ad}{ad}
\newcommand{\n}{\ensuremath{n}\xspace}
\renewcommand{\k}{\ensuremath{k}\xspace}

\newcommand{\Gc}[1]{\ensuremath{G^{\CC}_{#1}}}  
\newcommand{\Vc}[1]{\ensuremath{V^{\CC}_{#1}}}  
\newcommand{\gkn}{\Gc{\k,\n}\xspace}
\newcommand{\vkn}{\Vc{\k,\n}\xspace}
\newcommand{\kbein}{\ensuremath{\textmatrix{\Unity\\\Zero}}}
\newcommand{\Un}{\ensuremath{U(\n)}\xspace}
\newcommand{\Uk}{\ensuremath{U(\k)}\xspace}
\newcommand{\un}{\ensuremath{\mathfrak{u(\n)}}\xspace}

\newcommand{\lieu}{\ensuremath{\mathfrak{u}}}
\newcommand{\lieh}{\ensuremath{\mathfrak{h}}}

\renewcommand{\d}{\ensuremath{{d}}\xspace}
\newcommand{\dg}{\ensuremath{{d}^G}\xspace}
\newcommand{\dv}{\ensuremath{{d}^V}\xspace}
\renewcommand{\r}{\ensuremath{{r}}\xspace}
\newcommand{\rg}{\ensuremath{{r}^G}\xspace}
\newcommand{\rv}{\ensuremath{{r}^V}\xspace}
\newcommand{\noverk}{\frac{\n}{\k}}
\newcommand{\tnoverk}{\tfrac{\n}{\k}}
\newcommand{\varrhokc}{\obar{\varrho}}
\newcommand{\varrhouc}{\ubar{\varrho}}
\newcommand{\Dnk}{\ensuremath{D_{\n,\k}}}
\title{Sphere packing bounds in the Grassmann and Stiefel manifolds}
\author{Oliver Henkel (henkel@hhi.fraunhofer.de)\\
        Fraunhofer German-Sino Lab for Mobile Communications -- MCI\\
        Einsteinufer 37, 10587 Berlin, Germany
\thanks{Journal reference: IEEE Trans. Inform. Theory {\bf 51}, no 10 (2005),
3445--3456.
Personal use of this material is permitted. However, permission to
reprint/republish this material for advertising or promotional purposes or for
creating new collective works for resale or redistribution to servers or lists,
or to reuse any copyrighted component of this work in other works must be
obtained from the IEEE}
} 
\begin{document}
\maketitle
\begin{abstract}
   Applying the Riemann geometric machinery of volume estimates in terms of
   curvature, bounds for the minimal distance of 
   packings/codes in the Grassmann and Stiefel manifolds will be derived and
   analyzed. In the context of space time block codes this leads to a
   monotonically increasing minimal distance lower bound as a function of
   the block length. This advocates large block lengths for the code design.
\end{abstract}
%
%
%
%
%
%
%
%
%
%
%
%
\section{Introduction}
This work is inspired by Barg and Nogin's paper \cite{bar.nog} for
asymptotic packing bounds in the Grassmann manifold, based on an
asymptotic expression for the volume of metric balls. The basic
estimates defining the bounds are given by the well known
Gilbert--Varshamov and Hamming (or sphere packing) inequalities: In a
compact manifold $M$ without boundary furnished with a topological metric
$d$, let us denote the volume of the metric ball of radius $\delta$ as
$\vol\,B_d(\delta)$ (this quantity is presupposed to be independent of its
center). Then for any given $d_0$ there exists a packing (or code) 
${\cal C}\subset M$ with the prescribed minimal distance $\d_0$ and cardinality 
$\abs{{\cal C}}$ such that 
\begin{gather}
   \label{e.gilbert-varshamov-bound}
  \frac{\vol\,M}{\abs{\mathcal{C}}} \le \vol B_d(d_0) 
  \quad\text{(Gilbert-Varshamov)}\\
  \intertext{while for any packing/code $\mathcal{C}\subset M$ with data
    $(d_0,\abs{\mathcal{C}})$} 
   \label{e.hamming-bound}
  \vol B_d\big(\tfrac{1}{2}d_0\big) \le \frac{\vol\,M}{\abs{\mathcal{C}}}   
  \quad\quad\quad\quad\text{(Hamming)}
\end{gather}
holds.

Taking for $M$ the complex Grassmann manifold \gkn of $\k$ dimensional
complex subspaces of $\CC^{\n}$, Barg and Nogin derived closed form
expressions 
\begin{equation}
  \vol\,B_d(\delta) =
  \begin{cases}
     \left(\sin\frac{\delta}{\sqrt{\k}}\right)^{2\n\k+\o(n)} &
     \text{(geodesic distance)}\\
     \left(\frac{\delta}{\sqrt{\k}}\right)^{2\n\k+\o(n)} & 
     \text{(chordal distance )}
  \end{cases}
\end{equation}
as $\n\tends\infty$, leading to 
\begin{gather}
   \sqrt{\k}\arcsin\left( \frac{1}{\sqrt{2^{R/\k}}} \right) 
   \lesssim d_0 \lesssim 
   2\sqrt{\k}\arcsin\left( \frac{1}{\sqrt{2^{R/\k}}}  \right) \\
   \label{e.bn-chord-bounds}
   \sqrt{\frac{\k}{2^{R/\k}}} 
   \lesssim d_0 \lesssim 
   \sqrt{2\k\left(1-\left(1-\frac{1}{2^{R/\k}}\right)^2\right)}
\end{gather}
for geodesic, respectively chordal distance (defined later on), 
whereas $R$ denotes the rate 
\begin{equation}
   R = \frac{1}{\n}\log_2 \abs{\mathcal{C}}
\end{equation}

Furthermore Han and Rosenthal \cite{han.ros} recently derived upper bounds on
the minimal distance (more general: on the diversity of space time codes)
for packings on the unitary group \Un. 

A general capacity and performance analysis of space time codes in Rayleigh
flat fading MIMO scenarios without channel state information at the transmitter
\cite{hoc.mar.1,hoc.mar.2,zhe.tse,tar.jaf.cal} revealed that the
appropriate coding spaces are indeed 
\begin{itemize}
\item the (scaled) complex Grassmann
   manifold \gkn (set of \k dimensional linear subspaces of $\CC^{\n}$),
   if the channel is unknown at the receiver
\item the (scaled) complex Stiefel
   manifold \vkn (set of \k orthonormal vectors in $\CC^{\n}$) if the channel
   is known at the receiver. 
\end{itemize}
Here \k corresponds to the
number of transmit antennas and \n to the block length of the codes and the
work in \cite{bar.nog} refers to \gkn as $\n\tends\infty$ while
\cite{han.ros} refers to \vkn as $\k=\n$.

The aim of this work is to close the gap between those two results by
deriving bounds on the minimal distance for codes/packings in \gkn, \vkn
for arbritrary $(\k,\n)$ (section \ref{s.bounds}): 
Applying the bounds (\ref{e.gilbert-varshamov-bound}),
(\ref{e.hamming-bound}) with equality, 
the main task is to solve the equation 
\begin{equation}
   B_d(\delta)=c
   \mathcomma 
   c\in\RR
\end{equation} 
for (minimal) distances $\delta$ in \gkn, \vkn, with respect to some appropriate
distance measure \d. To this end volume estimates for the
volume of (small) balls $B_d(\delta)$ induced by curvature bounds for \gkn
and \vkn come into play. Associated comparison spaces with constant
curvature and simple volume forms provide bounds for $B_d(\delta)$. In
particular the lower bound turns out to permit a 
simple closed form expression with respect to $(\k,\n)$. Its analysis
culminates in Theorem \ref{thm.r} for the geodesic minimal distance lower
bound and Corollary 
\ref{cor.collected-results} for the minimal distance $\tilde{d}_0$ of the
corresponding space time codes. Surprisingly it turns out, that the minimal
distance $\tilde{d}_0$ grows at least proportional to $\sqrt{\n}$, while
keeping the rate and the transmit power per time step constant.
That is, increasing the block length enhances the possible minimal distance,
thus in coding spaces with large block lengths there exists codes with
potentially better error performance than in 'small' coding spaces. Since
most of the space 
time coding research efforts in the literature deal with small dimensional
coding spaces such as \Uk (e.g. \cite{tar.jaf.cal}), future research in the
more general \gkn, \vkn promises performance gains. 

Apart from space time codes recent developments in the design of space
frequency codes \cite{boe.bor.pau}, \cite{boe.bor} also indicate that the
relevant coding spaces are subspaces of large dimensional Stiefel and
Grassmann manifolds. Thus the achieved results here may be of considerable
importance for space frequency code design.
 
This article proceeds as follows. Section \ref{s.diffgeo} deals with
notational conventions and basic definitions concerning the Stiefel and
Grassmann manifolds (the coding spaces for space time or space frequency
codes).
In section \ref{s.bounds} explicit bounds for the minimal distance will be
calculated and compared to results obtained elsewhere. Further analysis
on the lower bound will be performed in section \ref{s.analysis}, culminating 
in Theorem \ref{thm.r}. Its implications for the minimal distance in coding
theory will be pointed out in Corollary
\ref{cor.collected-results}. Finally section \ref{s.conclusions} gives a
summary of the results.  
\section{The complex Stiefel and Grassmann manifolds}\label{s.diffgeo}
The complex Stiefel and Grassmann manifolds together with their topological
metrics (coding distance function in the language of coding theory) considered in
this work constitute the focus of this 
section. For the analysis in later sections we also need some explicit
curvature computations and rigorous proofs, which can be found in the
appendices \ref{app.diffgeo} and \ref{app.drequivalence}.

Readers who are mainly interested in the results concerning
packings/coding and who are willing to accept the
(quite standard) differential geometric facts can read this section without
reference to the appendices, where further details can be found.

A survey of the geometry of the \emph{real} Stiefel and
Grassmann manifolds aimed at non-specialists can be found in
\cite{ede.ari.smi}\footnote{The complex case considered here is similar to
  the real case but in some places certain peculiarities of the complex
  structure come into play} 
and for an elementary introduction to differential geometric
concepts see e.g. \cite{kue}.  
\subsection{The Stiefel manifold \vkn}\label{ss.stiefel-manifold}
The (complex) Stiefel manifold 
\begin{equation}\label{e.defvkn} 
   \vkn \mdef \{ \Phi\in\CC^{\n \times \k} \,|\, \Phi^{\dagger}\Phi=\Unity \}
\end{equation}
($\Unity$ denotes the identity matrix) 
can be equipped with the structure of an \Un-normal homogeneous space,
which justifies the coset representation 
\begin{equation}\label{e.vknunhomogeneous}
   \vkn \cong \Un \left/ \textmatrix{
                       \Unity & \Zero\\
                       \Zero & U(\n-\k)
                    }\right.\mathcomma
   \Phi \cong \bar{\Phi}\kbein
\end{equation}
($\bar{\Phi}\in\Un$), in particular
\begin{equation}\label{e.dimcvkn}
   \dimc\vkn=\dimc\Un-\dimc U(\n-\k)=\k(\n-\tfrac{\k}{2})
\end{equation}
and 
\begin{equation}\label{e.vol-stiefel}\begin{split}
   \vol\,\vkn & = \vol\,\Un / \vol\,U(\n-\k) \\
              & = \negthickspace\negthickspace
                  \prod_{i=\n-\k+1}^{\n} 
                  \negthickspace\negthickspace\negthickspace
                  \bigabs{S^{2i-1}}
                = \negthickspace\negthickspace
                  \prod_{i=\n-\k+1}^{\n} \frac{2\pi^i}{(i-1)!} 
\end{split}\end{equation}
For \vkn as a Riemannian manifold the concept of geodesics and geodesic
distance can be applied to obtain a canonical distance measure \rv:
Denoting the tangent space
of the unitary group \Un by \lieu(\n) consisting of skew-Hermitian
\n-by-\n matrices, tangents of \vkn may be represented as 
\begin{equation}\label{e.Vhorizontal} 
   \lieu(\n)\ni
   X=\eqmatrix{A & -B^{\dagger}\\B & \Zero}
   \mathcomma A\in\lieu(\k),\, B\in\CC^{(\n-\k)\times\k}
\end{equation} 
and 
\begin{equation}\label{e.rv}
   (\rv)^2=\oneover{2}\fnorm{X}^2=\oneover{2}\fnorm{A}^2+\fnorm{B}^2
\end{equation}
is the squared geodesic length of the geodesic connecting
$\Psi=\kbein\in\vkn$ with $\Phi=(\exp X) \kbein\in\vkn$. Here $\exp$
denotes the matrix exponential and the geodesic distance between arbitrary
points $\Psi',\Phi'\in\vkn$ follows from the isometric transformation
$\Psi=\bar{\Psi}'^{-1}\Psi'$ and $\Phi=\bar{\Psi}'^{-1}\Phi'$.

The canonical embedding \eqref{e.defvkn} of \vkn into the vector space 
$\left(\CC^{\n \times \k},\scalprod{\cdot,\cdot}_{\CC}\right)$ motivates the
definition of another topological ('chordal') metric/distance
\begin{equation}
   \dv(\Phi,\Psi) \mdef
     \fnorm{\Phi-\Psi} \mathcomma \Phi,\Psi\in\vkn
\end{equation}
which is important for space time coding, where it represents the decision
criterion at the maximum-likelihood receiver, if the channel is known at the
receiver ('coherent' case), see \cite{hoc.mar.2,tar.jaf.cal}. 
Note, that $\dv$ is entirely different from the geodesic distance \rv.
Nevertheless we have\footnote{At first sight the proposition seems 
  obvious, but one has to take into account that \dv is expressed in terms
  of $\Phi,\Psi\in\vkn$, while \rv is expressed in terms of the space of
  tangents and these two spaces are linked by the matrix exponential which
  can not be written in closed form compare Appendix
  \ref{app.drequivalence}. Furthermore unlike \dv, \rv is NOT induced by
  (geodesics with respect to) the seemingly canonical embedding
  $\vkn\subset\CC^{\n\times\k}$, compare Appendix \ref{sss.suppstiefel}} 
\begin{proposition}\label{prop.localequivdgv}\textin
   For $\k=\n$ or $\k\le\frac{\n}{2}$ the metrics $\dv$ and $\rv$ are
   locally equivalent, thus in sufficiently small neighborhoods there exist 
   constants $\alpha^V>0$, $\beta^V>0$ such that
   \begin{equation}\label{e.vchord-vs-geodesic}
      \beta^V \dv \le \rv \le \alpha^V \dv
   \end{equation}
   holds
\end{proposition}
This equivalence links the abstract (geodesic) sphere packing problem to
space time coding. The restriction to the cases $\k=\n$ and
$\k\le\frac{\n}{2}$ is mainly for convenience, since the main analysis will
concentrate on $\k\ll\n$.

\begin{proof}
   Lemma \ref{lem.dler}, \ref{lem.dgerkeqn}, \ref{lem.dger}
   in Appendix \ref{app.drequivalence}
\end{proof}
\begin{remark}\textin
   While it is an easy exercise to find 
   $\beta^V=\oneover{\sqrt{2}}$ (Lemma \ref{lem.dler}), 
   no concrete values for $\alpha^V$ have been obtained 
   rigorously. However, for $\k=2$, $\n=4,6,8$ numerical simulations led to 
   $\alpha^V\approx\frac{\pi}{2\cdot 0.9}$.
\end{remark}
\subsection{The Grassmann manifold \gkn}\label{ss.grassmann-manifold}
The (complex) Grassmann manifold 
\begin{equation}\label{e.defgkn}
   \gkn \mdef \{ \generate{\Phi} \,|\, \Phi\in\vkn \}
\end{equation}
of all $\k$-dimensional linear subspaces $\generate{\Phi}$ of
$\CC^{\n}$ also carries the structure of a \Un-normal homogeneous space
with coset representation
\begin{equation}\label{e.gknunhomogenous}
   \gkn \cong \Un \left/ \textmatrix{
                           \Uk   & \Zero\\
                           \Zero & U(\n-\k)
                         }\right.\mathcomma
   \generate{\Phi} \cong \Phi{\Phi^1}^{-1}
\end{equation}
(with $\Phi^1 \mdef (\Unity,\Zero)\Phi$) and
\begin{equation}\label{e.dimcgkn}
   \dimc\gkn=\k(\n-\k)
\end{equation}
The total volume of \gkn is
\begin{equation}\label{e.vol-grass}\begin{split}
   \vol\,\gkn & = \vol\,\vkn / \vol\,\Uk \\
              & = \negthickspace\negthickspace
                  \prod_{i=\n-\k+1}^{\n} \frac{2\pi^i}{(i-1)!} \left/
                  \prod_{j=1}^{\k} \frac{2\pi^j}{(j-1)!} \right.
\end{split}\end{equation}

Tangents become
\begin{equation}\label{e.Ghorizontal}
   X=\eqmatrix{\Zero & -B^{\dagger}\\B & \Zero}
   \mathcomma B\in\CC^{(\n-\k)\times\k}
\end{equation} 
with squared geodesic length 
$\oneover{2}\fnorm{X}^2=\fnorm{B}^2$, but 
there is an alternative notation in terms of the vector of principal angles
$\vartheta$ between subspaces:
To simplify matters let us assume $\k\le\n/2$ whenever we are in contact
with the Grassmann manifold. This is 
no restriction, since for $\k\ge\n/2$ we can always switch to the orthogonal
complement. Then there are precisely $\k$ principal angles 
$\vartheta_i$ between the subspaces
$\generate{\kbein}$ and $\generate{(\exp X)\kbein}$.
Performing a singular value decomposition on the  tangents
\eqref{e.Ghorizontal} one obtains (compare \ref{sss.suppgrassmann}).
$\fnorm{B}=\lnorm{\vartheta}{2}$ 
thus the geodesic distance \rg between $\generate{\kbein}$ and 
$\generate{(\exp X)\kbein}$ reads
\begin{equation}\label{e.rg}
   \rg=\oneover{\sqrt{2}}\fnorm{X}=\fnorm{B}=\lnorm{\vartheta}{2}
\end{equation}
As for \vkn there is also a different distance measure \dg in \gkn induced by
the maximum-likelihood receiver, which can be derived from the following 
geometric picture:

\vspace{1ex}{\em Spherical embedding:}
Unlike for the Stiefel manifold, there is no canonical embedding of \gkn
into Euclidean space unless choosing a representing unitary frame $\Phi_0$
in each subspace $\generate{\Phi}\in\gkn$. Nevertheless
there exists an interesting embedding of \gkn into Euclidean space given in
\cite{con.har.slo}: 
For $\Phi\in\vkn$ there is an well-defined associated orthogonal
projection  
\begin{equation}\label{e.PPhi}
   P_{\Phi}:=\Phi\Phi^{\dagger}:\mapset{\CC^{\n}}{\generate{\Phi}}
\end{equation}
of norm $\fnorm{P-\k/\n\,\Unity}^2=\k(\n-\k)/\n$ and $\tr P=\k$, which
justifies the embedding 
\begin{equation}\label{e.gspherical-embedding}\begin{gathered}
   \gkn \hookrightarrow S^{\n^2-2}\bigl(\sqrt{\k(\n-\k)/\n} \bigr)
      \subset\RR^{\n^2-1} \,, \\
   \mapto{\generate{\Phi}}{P_{\Phi}-\tfrac{k}{n}\,\Unity}
\end{gathered}\end{equation}
This motivates the 'chordal' topological metric
\begin{equation}
   \dg(\generate{\Phi},\generate{\Psi})\mdef
   \lnorm{\sin\vartheta}{2}=
   \frac{1}{\sqrt{2}}\fnorm{P_{\Phi}-P_{\Psi}}
\end{equation}
($\Phi,\Psi\in\vkn$).
Comparing \dg with the geodesic distance $\rg$ \eqref{e.rg}
between two subspaces we observe
\begin{proposition}\textin
   \begin{equation}\label{e.gchord-vs-geodesic}
      \beta^G \dg \le \rg \le \alpha^G \dg
   \end{equation}
whereas 
$\beta^G=1$ and
$\alpha^G = \frac{\pi}{2}$.
\end{proposition}
\section{Bounds for the minimal distance}\label{s.bounds}
Now let us specialize the general packing/coding bounds
\eqref{e.gilbert-varshamov-bound},\eqref{e.hamming-bound}. Set
\begin{equation}
   (M,\d) \mdef
   \begin{cases}
      (\vkn,\dv)\\
      (\gkn,\dg)
   \end{cases}
\end{equation}
and (compare \eqref{e.dimcvkn}, \eqref{e.dimcgkn})
\begin{equation}\label{e.dimensionD}
   D \mdef \dimr M = 
     \begin{cases} 
        \k(2\n-\k) \mathcomma M=\vkn\\
        2\k(\n-\k) \mathcomma M=\gkn
     \end{cases}
\end{equation}
for the two cases of interest. In the sequel other symbols like 
$\alpha$ are used generically to denote $\alpha^V$ or $\alpha^G$ when
specialized to the corresponding spaces \vkn, \gkn.
Denote by
\begin{equation}
   v(r) \mdef \vol B(r)
\end{equation}
the volume of the geodesic ball of radius $r$ in $M$, which is independent of its
center by left invariance of the Riemannian metric. With this notation the
Gilbert-Varshamov 
\eqref{e.gilbert-varshamov-bound} and Hamming bound 
\eqref{e.hamming-bound} for packings ${\cal C}$ on $M$ can
be compactly rewritten as
\begin{subequations}\label{e.specialized-bounds}
\begin{equation}\label{e.specialized-bounds-r}
     \ubar{\r}_0 \mdef v^{-1}\Big(\frac{\vol\,M}{2^{\n R}}\Big) 
     \le \r_0 \le 2 
     v^{-1}\Big(\frac{\vol\,M}{2^{\n R}}\Big) \mfed \bar{\r}_0
\end{equation}
\text{or relaxed w.r.t. coding distances $\d_0$ using
  \eqref{e.vchord-vs-geodesic}, \eqref{e.gchord-vs-geodesic}}
\begin{equation}\label{e.specialized-bounds-d}
     \ubar{\d}_0 \mdef \frac{1}{\alpha} v^{-1}\Big(\frac{\vol\,M}{2^{\n R}}\Big) 
     \le \d_0 \le
     \frac{2}{\beta} v^{-1}\Big(\frac{\vol\,M}{2^{\n R}}\Big) \mfed \bar{\d}_0
\end{equation}
\end{subequations}
So packing bounds are related to the coding bounds by simply setting
$\alpha=\beta=1$, thus replacing the (topological) metric distances by
geodesic distances. Due to the rather difficult to obtain explicit value
for $\alpha^V$ in \eqref{e.vchord-vs-geodesic} we focus on the packing
bounds \eqref{e.specialized-bounds-r} for most of the remaining analysis,
keeping in mind the simple relationship between statements about packings and
statements about space time coding.

To obtain the desired bounds for the minimal distance provided by
\eqref{e.specialized-bounds} we need closed form expressions for the volume
$v$ of small balls in $M$. As has been already indicated in the
introduction, this is a difficult task in general:
The canonical volume forms on \gkn and \vkn are elaborate to calculate.
Alternatively a common tool to compute volumes 
in Riemannian geometry arises from curvature, using Jacobi vector fields (see
e.g. \cite{gal.hul.laf} for details). Unfortunately a direct application 
can not be performed since we would have needed a diagonalization of
$XY-YX$ for each horizontal (compare \ref{ss.unitary})
$\norm{X}=\norm{Y}=1$ in $\lieu(\n)$ written in closed form.  
But there are simple volume estimates which will be presented in
\ref{ss.coarse-bounds}. In \ref{ss.comparisons} the results will be
compared to those already obtained 
in \cite{bar.nog}, \cite{han.ros}, in a few (computational simple) cases.
\subsection{Bishop/G\"unther volume bounds}\label{ss.coarse-bounds}
The method for volume computations in Riemannian manifolds using Jacobi
vector fields can be looked up in \cite[theorem 3.101]{gal.hul.laf}.
For $\kappa\in\RR$ let
\begin{equation}\label{e.volkappa}
   v^{\kappa}(r) \mdef
   \left(\frac{1}{\sqrt{\kappa}}\right)^{D-1}\bigabs{S^{D-1}}
   \int_0^r (\sin\sqrt{\kappa}t)^{D-1}\; dt
\end{equation}
denote the volume of the geodesic ball of radius $r$ in the manifold of
constant curvature $\kappa$ and let  
$\ubar{\kappa}\le\kappa\le\obar{\kappa}$ be defined by (compare
\eqref{e.ricci-curvature}, \eqref{e.sectional-curvature})
\begin{equation}\begin{gathered}
   \ubar{\kappa}\mdef\frac{1}{D-1}
     \min_{e_i}\ric(e_i,e_i) \\
   \obar{\kappa}\mdef
     \max_{\norm{X}=\norm{Y}=1} K(X,Y)
\end{gathered}\end{equation}
then we obtain monotone
volume bounds $v_l(r)\le v(r) \le v_u(r)$ for arbitrary $0\le r \le
\frac{\pi}{\sqrt{\rule{0pt}{1ex}\obar{\kappa}}}$ by 
\begin{equation}\label{e.vlu}
   v_l(r) \mdef v^{\obar{\kappa}}(r)
   \mathcomma\quad
   v_u(r) \mdef v^{\ubar{\kappa}}(r)
\end{equation}
From
$\lambda\le\kappa\followshort v^{\kappa}(r)\le v^{\lambda}(r)$ and 
$K\ge 0$ in $M$ \eqref{e.Kge0} we can further
relax $\ubar{\kappa}$ to zero, which yields the simple upper volume bound
\begin{equation}\label{e.vol-flatball-r}
   v_u(r) = v^0(r) = \bigabs{B^D} r^D
\end{equation}
A lower volume bound comes from an upper bound $\obar{\kappa}$ for $K$. 
Inserting  tangents $X, Y$ \eqref{e.Vhorizontal},
(resp. \eqref{e.Ghorizontal}) into (\ref{e.sectional-curvature}) subject to 
$\norm{X}=\norm{Y}=1$ yields
\begin{equation}
   K(X,Y) \le \obar{\kappa} =
     \begin{cases}
        2 & (\Uk=\Vc{\k,\k}) \\
        \frac{5}{2} & (\vkn, \k<\n)\\
        4 & (\gkn)
     \end{cases}
\end{equation}
Plugging this bounds into \eqref{e.specialized-bounds-r} we we end up with
\begin{equation}\label{e.coarsebounds}
     \ubar{\r}_0= (v^0)^{-1}\Big(\frac{\vol\,M}{2^{\n R}}\Big)
     \le \r_0 \le
     2 (v^{\obar{\kappa}})^{-1}\Big(\frac{\vol\,M}{2^{\n R}}\Big)=\bar{\r}_0
\end{equation}
With these settings an explicit (Maple-) calculation revealed
{\samepage
  {\par\footnotesize\[
    \arraycolsep0mm
    \begin{array}{c|cccc}
 k\backslash n & k & 2k & 3k & 4k\\\hline
1& [1.57,3.14]& [1.06,3.49]& [0.941,-1]& [0.886,-1]\\
 &            & \{0.500,1.05\}& \{0.595,1.40\}& \{0.630,1.71\}\\
\hline
2& [1.58,-1]& [1.38,-1]& [1.32,-1]& [1.29,-1]\\
 &          & \{0.771,-1\}& \{0.909,-1\}& \{0.973,-1\}\\
\hline
3& [1.74,-1]& [1.66,-1]& [1.63,-1]& [1.61,-1]\\
 &          & \{0.977,-1\}& \{1.15,-1\}& \{1.24,-1\}\\
\hline
4& [1.92,-1]& [1.92,-1]& [1.89,-1]& [1.88,-1]\\
 &          & \{1.15,-1\}& \{1.36,-1\}& \{1.46,-1\}
\end{array} 

    \]}
  \centerline{\footnotesize $[\ubar{\r}_0,\bar{\r}_0]$ for \vkn and
  $\{\ubar{\r}_0,\bar{\r}_0\}$ for \gkn with respect to
  \eqref{e.coarsebounds} for $R=1$}
}
{\samepage
  {\par\footnotesize\[
    \arraycolsep0mm
    \begin{array}{c|cccc}
 k\backslash n & k & 2k & 3k & 4k\\\hline
1& [0.0031,0.0061]& [0.0165,0.0330]& [0.0223,0.0446]& [0.0251,0.0502]\\
 &                & \{0.001,0.002\}& \{0.0055,0.0110\}& \{0.0098,0.0197\}\\
\hline
2& [0.0700,0.140]& [0.172,0.348]& [0.203,0.412]& [0.217,0.441]\\
 &          & \{0.0341,0.0682\}& \{0.0877,0.176\}& \{0.122,0.245\}\\
\hline
3& [0.217,0.440]& [0.416,0.898]& [0.467,1.04]& [0.490,1.11]\\
 &          & \{0.122,0.246\}& \{0.242,0.504\}& \{0.309,0.664\}\\
\hline
4& [0.403,0.847]& [0.678,-1]& [0.743,-1]& [0.771,-1]\\
 &              & \{0.242,0.504\}& \{0.422,0.992\}& \{0.517,1.75\}
\end{array} 

    \]}
  \centerline{\footnotesize $[\ubar{\r}_0,\bar{\r}_0]$ for \vkn and
  $\{\ubar{\r}_0,\bar{\r}_0\}$ for \gkn with respect to
  \eqref{e.coarsebounds} for $R=10$}
}\myline
($-1$ in the tables means, that Maple could not find a solution, due to
approximation error/too large sphere radii).
Observe that in the lower rate regime $\ubar{\r}_0$ still grows with $\n$
in \gkn, but slowly decreases in \vkn, while in the high rate regime
$\ubar{\r}_0$ is strict monotone with respect to $\n$, expecting the
intervals to become disjoint. 

So, while the high rate requirement is too restrictive, the results
for low rates are unsatisfactory in part. But the general
analysis of the lower bound in section \ref{s.analysis} will come up with
interesting results, supporting this approach. 
%
%
%
%
%
%
To clarify the presentation let us summarize the results so far in the
\begin{proposition}\label{prop.tabled-bounds}\textin
   The inequalities \eqref{e.coarsebounds}
   provide approximate bounds on the (geodesic) minimal distance for
   packings/space time codes on the Stiefel (coherent case) and Grassmann
   (non-coherent case) manifolds for any admissible $(\k,\n)$. 
   In particular the lower bound in \eqref{e.coarsebounds} is computational
   simple and guarantees the existence of corresponding packings/codes.
\end{proposition}
Especially for the Stiefel manifold these explicitly calculated bounds
appear to be new in the context of coherent space time coding.
\subsection{Comparisons with related results in the literature}
\label{ss.comparisons}
Han and Rosenthal \cite{han.ros} obtained bounds on the scaled chordal
distance $\Delta:=\frac{\dv}{2\sqrt{\k}}$ in the unitary
case $\k=\n$, $\Vc{\k,\k}=\Uk$. Based on a
numerically calculated exact volume they extracted
three upper bounds. The following table shows their (best) upper bounds
(2nd row) for $\Delta={\dv}/(2\sqrt{2})$ in $U(2)$
(in part relying on the results in \cite{lia.xia}) for
different rates $R$ (1st row) together with the upper bounds obtained here
(3rd row)

{\par\footnotesize\[
 \begin{array}{c|cccccccc}  
  R &                      
    2.29& 2.79& 3.0& 3.16& 3.32& 3.45& 3.50& 4.98\\\hline\noalign{\medskip} 
  \!\!\!\!\Delta\!\! & 
    0.675& 0.619& 0.597& 0.580& 0.558& 0.542& 0.535& 0.327\\\hline\noalign{\medskip} 
  \tfrac{\rv}{2} & 
    1.40& 1.01& 0.909& 0.843& 0.785& 0.742& 0.727& 0.409
\end{array}

\]}

Note, that equality in the (rough) estimate
$\frac{\dv}{2\sqrt{2}}\le\frac{\rv}{2}$ (Lemma \ref{lem.dler}) has been
forced in the third row of the table to convert the geodesic distances
computed by a Maple program into chordal distances.
Consequently the bounds of \cite{han.ros} are tighter than the bounds
obtained here, since in the case of unitary matrices there are 
more specialized (but less general) methods available to obtain bounds.

As already stated in the introduction another (asymptotic) result has been
obtained by Barg and Nogin. For the non-asymptotic case they presented an
exact volume formula \cite[eq. (11)]{bar.nog} for regions in the 
(real and complex) Grassmann manifold. 
\begin{equation}
   \label{e.concrete-exact-gball-volume}
   \begin{split}\raisetag{2cm}
   \vol B(r) 
      & = 2^{\k}\abs{\gkn}
          \prod_{i=1}^{\k}\frac{(\n-i)!}{[(i-1)!]^2 (\n-\k-i)!} \times \\
      &   \int_{\substack{0<\vartheta_1<\dots<\vartheta_{\k}<\pi/2\\
                          \norm{\vartheta}_2\le r}}\,
          d\vartheta_1\dots d\vartheta_{\k} \times \\
      &   \prod_{i=1}^{\k}(\sin\vartheta_i)^{2(\n-2\k)+1}\cos\vartheta_i\,
          \prod_{j<l}(\sin^2\vartheta_j-\sin^2\vartheta_l)^2
\end{split}\end{equation}
which can be computed in polar coordinates, compare \ref{appss.volumes}.
Although \eqref{e.concrete-exact-gball-volume} is exact, it does not provide a
closed form for varying dimensions. Moreover the computations are elaborate
compared with the ones done here, such that the evaluation of
(\ref{e.specialized-bounds}) become intractable. 
%
%
%
%
%
%
%
%
%
%
%
\section{Analysis of the lower bound}\label{s.analysis}
The lower bound for the (geodesic) minimal distance $\r_0$ guarantees the
existence of corresponding packings/codes. Due to \eqref{e.vol-flatball-r}
we can explicitly solve the lower bound in \eqref{e.coarsebounds}:
\begin{equation}
   \ubar{\r}_0^{\Dnk} = \oneover{2^{\n R}} \cdot
              \frac{\vol\,M}{\bigabs{B^{\Dnk}}}          
\end{equation}
with \Dnk\ defined as $D$ in \eqref{e.dimensionD}, thus
\begin{equation}
   \label{e.Dnk}\tag{\ref{e.dimensionD}'}
   \Dnk = 2\n\k-\epsilon\k^2
   \mathcomma
   \epsilon=\begin{cases} 1,\,\vkn \\ 2,\,\gkn \end{cases}
\end{equation}
Then
{\samepage
\begin{theorem}\label{thm.r}\textin
   The (geodesic) minimal distance $\r_0$ in $M$ can be lower bounded by
   \begin{equation}\label{e.r0-lower}
      \ubar{\r}_0 
         \ge \left( \oneover{2} \right)^{\frac{\n R}{D_{\n,\k}}}
   \end{equation}
   with the right hand side monotonically increasing as a function of
   $\n$ for $\n\ge\k$. Asymptotically 
   \begin{equation}\label{e.r0-limit}
      \lim_{\n\tendshort\infty}\ubar{\r}_0 
        = \sqrt{\frac{\k}{2^{R/\k}}}
   \end{equation}\vspace{-2ex}
   holds.
\end{theorem}
\textout
}
In particular this establishes a \emph{monotonically increasing} lower
estimate for $\ubar{\r}_0$ common for \vkn and \gkn, which is
not obvious from the picture drawn 
from the explicit calculations of $\ubar{\r}_0$ for rate $R=1$ in the
previous sections.
Of course, the theorem also holds for the (topological) minimal distance 
$\ubar{\d}_0=\oneover{\alpha}\ubar{\r}_0$, connecting this result with
space time coding theory.

\textout
\begin{proof}\myline
   Set 
   $a\mdef 2^{-\n R/\Dnk}$ and
   $b\mdef\big(\frac{\vol\,M}{\abs{B^{\Dnk}}}\big)^{1/\Dnk}$.
   Then $\ubar{\r}_0 = a\cdot b$ and from \eqref{e.Dnk}
   $a$ is monotonically increasing as a function of \n with
   $\lim_{\n\tendshort\infty} a = 2^{-R/2\k}$.\myline
   For $b$ we show
   \begin{equation}\label{e.b-estimate}
      b \ge 1
   \end{equation}
   and
   \begin{equation}\label{e.b-limit}
      \lim_{\n\tendshort\infty} b = \sqrt{\k}
   \end{equation}
   and the theorem follows.\\[1ex]
   For the two cases of interest $b^{\Dnk}$ is given as (using
   \eqref{e.vol-flatball},\eqref{e.vol-stiefel},\eqref{e.vol-grass}) 
   \begin{equation}
      b^{\Dnk} = 
        \left\{
          \begin{aligned}
           &\frac{\vol\,\vkn}{\bigabs{B^{\k(2\n-\k)}}}
             = (2\sqrt{\pi})^{\k}\frac{\Gamma(\k(2\n-\k)/2+1)}
                                      {\prod_{i=\n-\k+1}^{\n}\Gamma(i)}
           \\
           &\frac{\vol\,\gkn}{\bigabs{B^{2\k(\n-\k)}}}
             = \frac{\Gamma(2\k(\n-\k)/2+1)\prod_{j=1}^{\k}\Gamma(j)}
                    {\prod_{i=\n-\k+1}^{\n}\Gamma(i)}
        \end{aligned}
        \right.
   \end{equation}
   for $\n\ge\k$, resp. $\n\ge\k+1$.
   The proof of \eqref{e.b-estimate} relies on the simple estimate
   \begin{equation}\label{e.gamma-estimate}\begin{split}
      \frac{\Gamma(M+1)}{\Gamma(m+1)}
      & = (m+1)(m+2)\cdot\,\dots\, \cdot M \\
      & \ge (m+1)^{M-m}
        \ge m^{M-m}
   \end{split}\end{equation}
   for $m,M\in\toneover{2}\NN$,$M-m\in\NN$.
   Since $\Dnk>0$ it suffices to show $B_{\k,\n}\mdef b^{\Dnk} \ge 1$ for all
   admissible $(\k,\n)$. This will be proven by induction over $\k$ and
   $\n$.\\[1ex]
   $\boldsymbol{\vkn}:$
      \begin{enumerate}
      \item $B_{1,1} = 2\sqrt{\pi}\,\Gamma(3/2)=\pi>1$
      \item Induction over $\k$
         \[\begin{split}
           \frac{B_{\k+1,\k+1}}{B_{\k,\k}} 
            &=2\sqrt{\pi}\frac{\Gamma(\overbrace{\k^2/2+(\k+1/2)}^{M}+1)}
                              {\Gamma(\k^2/2+1)\Gamma(\k+1)} \\
            &\underset{\eqref{e.gamma-estimate}}{\ge} 
              2\sqrt{\pi}\frac{(\k^2/2+1)^{\k+1/2}}{\k!}\\
            &\underset{\k^2/2+1>\k}{\ge}
              2\sqrt{\pi}\sqrt{\k^2/2+1}\\
            &\underset{\k=1}{\ge} \sqrt{6\pi} > 1
         \end{split}\]
         Thus $B_{\k,\k}>1\: \forall_{\k\ge 1}$
      \item Induction over $\n\ge\k$
         \[\begin{split}
           \frac{B_{\k,\n+1}}{B_{\k,\n}}
            &=\frac{\Gamma(\k/2(2\n-\k)+\k+1)}
                   {\Gamma(\underbrace{\k/2(2\n-\k)}_{m}+1)}
              \frac{\Gamma(\overbrace{\n-\k}^{\tilde{m}}+1)}{\Gamma(\n+1)}\\
            &\underset{\eqref{e.gamma-estimate}}{=}
              \frac{(m+1)\cdot\,\dots\,\cdot(m+\k)}
                   {(\tilde{m}+1)\cdot\,\dots\,\cdot(\tilde{m}+\k)}
             \ge 1 
         \end{split}\]
         since 
         $m-\tilde{m}=\n(\k-1)-\k^2/2+k
          \underset{\n\ge\k}{\ge}\k^2/2>0$.
         Thus for every $\k\ge1$ we have $B_{\k,\n}>1\:\forall_{\n\ge\k}$.
      \end{enumerate}
   $\boldsymbol{\gkn}:$
      \begin{enumerate}
      \item $B_{1,2} = \frac{\Gamma(2)\Gamma(1)}{\Gamma(2)} = 1$
      \item Induction over $\k$
         \[\begin{split}
           \frac{B_{\k+1,\k+2}}{B_{\k,\k+1}} 
             =\frac{\Gamma(\k+2)}{\Gamma(\k+1)}
              \Gamma(\k+1)
              \frac{\prod_{i=2}^{\k+1}\Gamma(i)}
                   {\prod_{i=2}^{\k+2}\Gamma(i)}
             = 1
         \end{split}\]
         Thus $B_{\k,\k+1}=1\: \forall_{\k\ge 1}$
      \item Induction over $\n\ge\k+1$\vspace{-2ex}
         {\samepage
         \[\begin{split}
           \frac{B_{\k,\n+1}}{B_{\k,\n}}
            &=\frac{\Gamma(\k(\n-\k)+\k+1)}
                   {\Gamma(\underbrace{\k(\n-\k)}_{m}+1)}
              \frac{\Gamma(\overbrace{\n-\k}^{\tilde{m}}+1)}{\Gamma(\n+1)}\\
            &\underset{\eqref{e.gamma-estimate}}{=}
              \frac{(m+1)\cdot\,\dots\,\cdot(m+\k)}
                   {(\tilde{m}+1)\cdot\,\dots\,\cdot(\tilde{m}+\k)}
             \ge 1
         \end{split}\]
         since 
         $m-\tilde{m}=\n(\k-1)-\k^2+k
          \underset{\n\ge\k+1}{\ge}\k-1\ge0$}
         and it follows for every $k\ge1$, that
         $B_{\k,\n}\ge1\:\forall_{\n\ge\k+1}$ as desired. 
      \end{enumerate}
   Let us now prove \eqref{e.b-limit}. At first 
   $(2\sqrt{\pi})^{\k/\Dnk}\underset{\n\tendshort\infty}{\tends} 1$ and
   $\prod_{j=1}^{\k}\Gamma(j)^{1/\Dnk}\underset{\n\tendshort\infty}{\tends} 1$
   holds. So it remains the evaluation of 
   \[
      \lim_{\n\tendshort\infty}
        \left(\frac{\Gamma(\Dnk/2+1)}{\prod_{i=\n-\k+1}^{\n}\Gamma(i)} 
        \right)^{1/\Dnk}
   \]
   Stirling's formula reads either ('$\sim$' denotes asymptotic equivalence)
   \begin{equation*}
      \Gamma(m+1) \sim \sqrt{2\pi m}\left( \frac{m}{e} \right)^m 
      \quad\text{or}\quad
      \Gamma(m) \sim \sqrt{\frac{2\pi}{m}}\left( \frac{m}{e} \right)^m 
   \end{equation*}
   and by $\Dnk\sim 2\n\k$ we deduce
   \[\begin{split}
     &\Big( \frac{\Gamma(\Dnk/2+1)}{\prod_{i=\n-\k+1}^{\n}\Gamma(i)} 
      \Big)^{1/\Dnk}\hspace{-2ex}
      \sim
        \left(
          \frac{\sqrt{2\pi \Dnk/2}\sqrt{\tfrac{\Dnk}{2e}}^{\Dnk}}
               {\prod_{i=\n-\k+1}^{\n}
                  \sqrt{\tfrac{2\pi}{i}}\big(\tfrac{i}{e}\big)^i}
        \right)^{1/\Dnk}\\
     &= \underbrace{\sqrt{\frac{\pi}{(2\pi)^{\k}}}^{1/\Dnk}}_{\tends 1}
        \oneover{\sqrt{2e}}
        \frac{\overbrace{\sqrt{\Dnk}^{1/\Dnk}}^{\tends 1} \sqrt{\Dnk}}
             {\left[\prod_{i=\n-\k+1}^{\n} 
                    \toneover{\sqrt{i}}\big(\tfrac{i}{e}\big)^i\right]^{1/\Dnk}}\\
     &\sim \oneover{\sqrt{2e}}
      \underbrace{\oneover{
          \left(\prod_{i=\n-\k+1}^{\n}\toneover{\sqrt{i}}\right)^{1/\Dnk}}}_
          {\tends 1}\; \times\\
     &\hspace{3ex}
      \underbrace{\sqrt{\Dnk}\,e^{-1/\Dnk\sum_{i=\n-\k+1}^{\n} i\ln\tfrac{i}{e}}}_
                 {\tends \sqrt{2e\k}} \\
     &\sim \sqrt{\k}
   \end{split}\]
   This proves \eqref{e.b-limit}
\end{proof}
\subsection{Final remarks and application to coding theory}
\label{ss.applications}
A remarkable coincidence arises from Barg/Nogin's results for the chordal
distance in \gkn. Denoting the lower bound in \eqref{e.bn-chord-bounds} by
$\ubar{\delta}_0$ we find
\begin{equation}
   \ubar{\delta}_0 
     = \lim_{\n\tendshort\infty}\ubar{\r}_0
\end{equation}
therefore, the geodesic lower bound $\ubar{\r}_0$ obtained from the flat
geodesic volume estimate $v^0(r)\le v(r)$ asymptotically equals the
(asymptotic) exact chordal lower bound \eqref{e.bn-chord-bounds}. This
seems reasonable since in flat space the geodesic distance coincides with the
Euclidean (chordal) distance.

Apart from the asymptotics, let us consider the lower bound
\eqref{e.r0-lower} of Theorem \ref{thm.r}. It guarantees the existence of
packings with minimal distance $\r_0$ bounded monotonically from
below in $\vkn$, resp. $\gkn$, when $\n$ grows. In coding
theory 
$\left(\widetilde{\vkn}=\sqrt{\noverk}\,\vkn\mathcomma\dv\right)$, resp.
$\left(\widetilde{\gkn}=\sqrt{\noverk}\,\gkn\mathcomma\dg\right)$,
represent the coding spaces
for space time block codes for the Rayleigh flat fading channel unknown to
the transmitter and known, resp. unknown, channel at the receiver. The factor
$\sqrt{\noverk}$ serves as a constraint, holding the transmit power at each
time step constant for different choices of $(\k,\n)$, thus provide a fair
comparison of codes from different coding spaces. 
In a Riemannian manifold $M$ with metric $g$ the mapping
$\mapto{(\lambda M,g)}{(M,\lambda^2 g)}$, $\lambda>0$, is isometric, leading 
immediately to the scaled geodesic minimal distance 
$\tilde{\r}_0=\lambda\r_0$. With respect to the coding
distances $\d$ we obtain instead
\begin{equation}\label{e.muhomothety}\begin{split}
 & \Big(\sqrt{\tnoverk} M,\d\Big) 
   \cong
   \Big(M,\sqrt{\mu\tnoverk}\,\d\Big)
   \mathcomma\\
 & \mu\mdef
   \begin{cases}
      1 \mathcomma\! M=\vkn\\
      \oneover{2}\mathcomma\! M=\gkn
   \end{cases}
\end{split}\end{equation}
whereas $\tfrac{\rho}{4}\mu\tnoverk$ ($\rho\ge 1$ denoting the signal to noise
ratio) is (a lower bound of) the first order term (the so called diversity
sum, our metric here) in the expansion 
of the Chernov bound for the pairwise error probability, compare
\cite[formulas (17)(18)(19)(20)]{hoc.mar.2}: 
The factor of $\oneover{2}$ for \gkn stems from the slightly different 'effective'
transmit power $\varrhouc\mdef\frac{(\rho\n/\k)^2}{4(1+\rho\n/\k)}$ compared to the
known channel effective transmit power $\varrhokc\mdef\frac{\rho\n}{4\k}$, satisfying 
$\oneover{2}\varrhokc \le \varrhouc \le \varrhokc$, whereas $\rho\ge
1$,$\n\ge\k$ is understood. 
Collecting all formulas we finally infer from Theorem \ref{thm.r}:
\begin{corollary}\label{cor.collected-results}\textin
Given $\rho\ge 1$ and $\n\ge 2\k$, there exist space time block codes with
minimal distance $\tilde{\d}_0$ lower bounded by 
\begin{equation}
   \tilde{\ubar{\d}}_0 
     = \frac{\sqrt{\mu}}{\alpha}\sqrt{\noverk} \ubar{\r}_0
     \ge \frac{\sqrt{\mu}}{\alpha}\sqrt{\noverk} 
         \left( \oneover{2} \right)^{\frac{\n R}{\Dnk}}
\end{equation}
whereas $\alpha$ is determined by \eqref{e.vchord-vs-geodesic},
resp. \eqref{e.gchord-vs-geodesic}, $\Dnk$ is defined in
\eqref{e.dimensionD} (resp. \eqref{e.Dnk}),
and $\mu$ in \eqref{e.muhomothety}.\myline
Thus the performance (which scales with ${\tilde{\d}_0}^2$) potentially
increases monotonically at least proportionally to $\noverk$.
\end{corollary}
The last statement in the corollary follows from the observation, that the
diversity (essentially the inverse of the Chernov bound for the pairwise
error probability) as a basic performance measure for space time codes
\cite{hoc.mar.2} is a homogenous polynomial. The first order term coincides
with the metric ${\tilde{\d}}^2$, while all higher order terms scale with
a power of ${\tilde{\d}}^2$ when code design is interpreted as a
constrained packing problem (considering the higher order terms as
constraints according to a normalized distance distribution).
\section{Conclusions}\label{s.conclusions}
The framework in \cite{bar.nog}, \cite{han.ros} had been successfully
generalized to the 
Stiefel manifold \vkn, $\n\ge\k$, and to \gkn, $\infty\gg\n\ge\k/2$ using the
completely different method of Riemannian volume bounds (Proposition
\ref{prop.tabled-bounds}). Unlike the exact volume formula the lower bound
can be relatively simple 
analyzed as a function of $(\k,\n)$ for both \gkn and \vkn, leading to
Theorem \ref{thm.r}, resp. Corollary \ref{cor.collected-results}. Although
the used estimates were quite conservative 
they apply (in principle) in any Riemannian homogeneous spaces. 

The connection to the coding theory of space time block codes
advocates further efforts in finding codes in the spaces \vkn, resp. \gkn
for $\n$ much larger than \k. Since the minimal distances grow
proportionally to $\sqrt{\noverk}$ while the transmit power per time step
remains constant, there is a considerable performance
impact to expect, when coding in \vkn (resp. \gkn) as opposed to coding in
$U(\k)$\footnote{
  Note that this does not contradict the (converse) conclusions in
  \cite{foz.mclau.schaf}, which do not apply here: The error probability
  computations done there with respect to increasing block length 
  $\n\tendshort\infty$ are constrained by a fixed total number of code
  symbols sent. This is a different scenario, not relevant for the analysis
  performed here.}.
Furthermore, as already pointed out in the introduction the
developments in space frequency coding indicate, that the relevant
coding spaces are subsets in some \vkn (resp. \gkn) whereas the number of
subcarriers $\n$ satisfies $\n\gg\k$, thus the results proven here may apply
to space frequency codes as well.

\section*{Acknowledgment}
I would like to thank Peter Jung for many helpful discussions.
\useRomanappendicesfalse
\appendices
\renewcommand{\theequation}{\thesection.\arabic{equation}}
\numberwithin{equation}{section}
\def\thesubsection{\thesection-\Roman{subsection}} 
\def\thesubsectiondis{\Roman{subsection}.} 
\section{Differential geometric calculations}\label{app.diffgeo}
\subsection{$\Un$-normal homogeneous spaces}
\label{appss.un-homogeneous-space}
For the theoretical background, common notation and curvature formulas
I refer to \cite{poor} as a reference.

\subsubsection{The unitary group}\label{ss.unitary}
\begin{equation}
   \Un = \{ \bar{\Phi} \in\CC^{\n\times\n} \,|\, 
            \bar{\Phi}^{\dagger}\bar{\Phi}=\Unity \}
\end{equation}
is a compact, connected Lie group and a real manifold of dimension
$\dimr=\n^2$. The corresponding Lie algebra (i.e. the tangent space of \Un
at $\Unity$) is
\begin{equation}
   \lieu(\n) = \{ \bar{X} \in \CC^{\n\times\n} \,|\,
                  \bar{X}^{\dagger} = -\bar{X} \}
\end{equation}
and the matrix exponential $\exp$ maps $\lieu(\n)$ into \Un.
On $\lieu(\n)$ the (bi-invariant) Riemannian metric for \Un is defined as
\begin{equation}
   \scalprod{\bar{X},\bar{Y}}=\frac{1}{2}\tr \bar{X}^{\dagger}\bar{Y}
\end{equation} 
thus $\scalprod{\bar{X},\bar{X}}=\oneover{2}\fnorm{\bar{X}}^2$
(Frobenius norm) holds. 

A manifold $M$ is a \Un-homogeneous space, if there is a
transitive \Un action on $M$ such that $M\cong\Un/H$ for some isotropy
subgroup $H\subset\Un$. If $\lieh\subset\lieu(\n)$ denotes the Lie algebra
of $H$ there is a canonical decomposition of tangent vectors 
$\lieh\oplus\lieh^{\smallperp}=\lieu(\n)
 \ni \bar{X}=X^{\smallparallel}+X$
and we can identify tangents of $M$ with so called 'horizontal' tangent vectors
$X\in\lieh^{\smallperp}$. With this identification $M$ is
called normal homogeneous. 
%

Then the sectional curvature $K$ and the Ricci curvature $\ric$ of $M$ are
given as  
\begin{gather}
   \label{e.sectional-curvature}
   K(X,Y) = \oneover{4}\norm{[X,Y]}^2
          + \frac{3}{4} \norm{[X,Y]^{\smallparallel}}^2\\
  \label{e.ricci-curvature}
   \ric(e_i,e_i) = \sum_j K(e_i,e_j)
\end{gather}
whereas $[X,Y]=XY-YX$, $X$ and $Y$ are normalized tangent vectors and
$\{e_i\}$ denotes a orthonormal base in $\lieh^{\smallperp}$.
Note that the sectional curvature $K$ is always non-negative
\begin{equation}\label{e.Kge0}
   K \ge 0
\end{equation}

\subsubsection{Supplements for the Stiefel manifold}\label{sss.suppstiefel}
The (complex) Stiefel manifold \eqref{e.defvkn}
is canonically a \Un-normal homogeneous space:  
The canonical left multiplication  
of $\k$-frames in $\CC^{\n}$ by unitary $\n \times \n$ matrices 
transforms each pair of $\k$-frames into each other. Thus the group action
of \Un on \vkn is transitive with isotropy group
$H=\textmatrix{\Unity & \Zero\\ \Zero & U(\n-\k)}$
and establishes the 
canonical diffeomorphism \eqref{e.vknunhomogeneous}.
Then $\lieh=\textmatrix{\Zero & \Zero\\ \Zero & \lieu(\n-\k)}$ and
 tangents $X\in\lieh^{\smallperp}$ have the form
\eqref{e.Vhorizontal}, and \eqref{e.rv} follows for the geodesic distance
\rv. Note that this distance is not induced by the
length of the geodesics obtained from the canonical embedding of \vkn into
$\CC^{\n\times\k}$, compare \cite{ede.ari.smi} in the real case and
additionally \cite[Example 6.61(b)]{poor} in the complex case.

\subsubsection{Supplements for the Grassmann manifold}\label{sss.suppgrassmann}
The (complex) Grassmann manifold \eqref{e.defgkn}
carries the structure of a \Un-normal homogeneous space by
forgetting 
not only the orthogonal complement of $\bar{\Phi}\in\Un$ (which has been
done for \vkn) but also 
the particular choice of the spanning $\k$-frame. 
Thus $H=\textmatrix{U(\k) & \Zero\\ \Zero & U(\n-\k)}$ and this leads to
\eqref{e.gknunhomogenous}. Note that the coordinate representation 
$\generate{\Phi} \cong \Phi{\Phi^1}^{-1}$ holds only locally in general,
but it turns out, that this representation covers all but a set of measure
zero, hence we abandon this distinction between local and global properties in
this work and drop the distinction between \gkn and its coordinate
domain. Calculating $\lieh^{\smallperp}$ leads to  tangents of the form
\eqref{e.Ghorizontal}.

Given two elements $\generate{\Phi},\generate{\Psi}\in\gkn$ then the
$\k$ stationary angles $0\le\vartheta_1\le\dots\le\vartheta_{k}\le\pi/2$ between 
$\generate{\Phi}$ and $\generate{\Psi}$ are defined
successively by the critical values
$\arccos\abs{\scalprod{v_i, w_i}}$, $i=1,\dots,\k$ (in increasing order),
of $\mapto{(v,w)}{\arccos\abs{\scalprod{v, w}}}$
where the unit vectors $v,w$ vary over 
$\{v_1,\dots,v_{i-1}\}^{\perp} \subset \generate{\Phi}$, 
respectively
$\{w_1,\dots,w_{i-1}\}^{\perp} \subset \generate{\Psi}$.
It is well known that the stationary angles 
can be computed by the formula (any representing $\k$-frame will do)
\begin{equation}\label{e.principalangles}
   \cos\vartheta_i = \sigma_i(\Phi^{\dagger}\Psi)
   \mathcomma
   i=1,\dots,\k
\end{equation}
whereas $\sigma_i(M)$, $i=1,\dots,\k$ denotes the $i$-th singular value of
the matrix $M$ in decreasing order.

Given a tangent $X=\textmatrix{\Zero & -B^{\dagger}\\ B & \Zero}$,
$B\in\CC^{(\n-\k)\times\k}$,  
the singular value decomposition 
$B=V\Sigma W^{\dagger}=V_1 S W^{\dagger}$, 
$V=(V_1,V_2)\in U(\n-\k)$, $W\in\Uk$,
$\Sigma=\textmatrix{S\\\Zero}$, $S=\diag(\sigma_1,\dots,\sigma_{\k})$,
yields $X=U \Delta U^{\dagger}$ with
$U=\textmatrix{W & 0 & 0 \\ 0 & V_1 & V_2}$, 
$\Delta=\textmatrix{{\cal D} & 0\\ 0 & 0}$, and
${\cal D}=\textmatrix{0 & -S\\ S & 0}$.
From this one calculates 
$(\exp X)\kbein=\textmatrix{W(\cos S)W^{\dagger} \\ 
                              V_1(\sin S)W^{\dagger}}$,
and 
$\cos\vartheta=\sigma\big(\textmatrix{\Unity & \Zero}(\exp X)\kbein\big)
 = \cos S$, 
thus $\vartheta=\sigma$
and \eqref{e.rg} follows.

The space of orthogonal projections 
$\Pi_V\mdef\{P_{\Phi} \,|\, \Phi\in\vkn\}$ (compare \eqref{e.PPhi}) can be
identified with \gkn. In particular we have $\Pi_V=\Pi_{\k}$ with
\begin{equation}\label{e.projection-operators}\begin{split}
   \Pi_{\k}\mdef \big\{P\in \CC^{\n \times \n} \,|\, & 
   P^{\dagger}=P, P^2=P, \tr P=\k, \\
   & \fnorm{P-\k/\n\,\Unity}^2=\k(\n-\k)/\n \big\}
\end{split}\end{equation}
as one can see by picking an appropriate representative
$\Phi\in\generate{\Phi}$ (e.g. $\Phi=\kbein$ due to invariance of
$\Pi_{\k}$ under the left and right unitary action).
Since each $P\in\Pi_{\k}$ is Hermitian with constant trace, $\Pi_{\k}$ is
canonically a real submanifold of $\RR^{\n^2-1}$, the constant norm
justifies the embedding \eqref{e.gspherical-embedding}
\subsection{Volume computations}\label{appss.volumes}
\subsubsection{Total volume}
The unitary group $\Un\subset\Gl(\n,\CC)$ can be equipped with the induced
Lebesgue measure from 
the ambient space $\RR^{2\n^2}$. The Stiefel manifold inherits its volume
measure from its total space \Un: 
We get from the familiar volume formulas 
\begin{gather}
   \label{e.vol-flatsphere}
   \bigabs{S^{m-1}}\mdef\vol\,S^{m-1} = \frac{2\pi^{m/2}}{\Gamma(m/2)}\\
   \label{e.vol-flatball}
   \bigabs{B^m}\mdef\vol\,B^m = \bigabs{S^{m-1}}/m
\end{gather}
for the unit sphere $S^{m-1}$ and the unit ball $B^m$ in $\RR^m$
and the canonical
homogeneous family $S^{2m-1} \cong U(m)/U(m-1)$
the following recursive formula 
$\vol U(1)=\abs{S^1}=2\pi$, 
$\vol U(m)=\abs{S^{2m-1}}\vol U(m-1)$, 
and therefore \eqref{e.vol-stiefel} and \eqref{e.vol-grass}.
\subsubsection{Volume for regions in \gkn}
The volume formula for regions in the complex Grassmann manifold will be
derived, based on \cite[eq. (11)]{bar.nog}\footnote{Unfortunately (in
  their first paper version) their formula is not 
  correct in the complex case (private communication).
  Fortunately this does not affect the (asymptotic) results obtained in
  \cite{bar.nog}. An erratum has already been produced, thus the derivation
  here is only for completeness of the presentation and the convenience of
  the reader}. 

Starting with formula \cite[(A.18)]{has.hoc.mar} for the distribution of
eigenvalues $\lambda_i$, with $\lambda_i=\cos^2\vartheta_i$ of 
$\textmatrix{\Unity,\Zero}\Phi\Phi^{\dagger}\kbein$ we obtain the volume
density as the marginal density 
\begin{equation}\label{e.volform}\begin{split}
   \omega & = C(\k,\n)\cdot\,
              \prod_{i=1}^{\k}(1-\lambda_i)^{\n-2\k}
              \prod_{j<l}(\lambda_l-\lambda_j)^2\cdot\,
              d\lambda_1\dots d\lambda_{\k} \\
          & = C(\k,\n)\cdot 2^{\k}k!\cdot\,
              \prod_{i=1}^{\k}(\sin\vartheta_i)^{2(\n-2\k)+1}\cos\vartheta_i \\
          &   \hspace*{2.75cm}\prod_{j<l}(\sin^2\vartheta_j-\sin^2\vartheta_l)^2\cdot\,
              d\vartheta_1\dots d\vartheta_{\k}
\end{split}\end{equation}
whereas the Jacobi determinant 
$2^{\k}\prod_{i=1}^{\k}\sin\vartheta_i\cos\vartheta_i$ 
of the mapping $\mapto{\lambda}{\vartheta}$ has been introduced in order to
express the volume density in terms of 
$\vartheta$, and $\k!$ establishes the ordering condition on the (open) simplex 
$\Theta=\{0<\vartheta_1<\cdots<\vartheta_{\k}<\pi/2\}$
of stationary angles. 
The constant $C$ is just a normalization factor, which reads in our case
\begin{equation}
   C(\k,\n)=\frac{\abs{\gkn}}{\k!}
            \prod_{i=1}^{\k}\frac{(\n-i)!}{[(i-1)!]^2 (\n-\k-i)!}
\end{equation}
(without the factor $\abs{\gkn}$ this would give the Haar measure used in
\cite{bar.nog} on \gkn).
The volume of sufficiently small geodesic balls is now given as
\begin{equation}\label{e.exact-gball-volume}\begin{split}
   & \vol B(r) 
          = \int_{\Theta\cap\{\norm{\vartheta}_2\le r\}} \omega(\vartheta) \\
        & = \int_0^r d\rho 
            \int_{\alpha,\beta_i\in[0,\pi/2]}
                 \oneover{\k!}\omega(\vartheta(\rho,\alpha,\beta))
                 \abs{\det J_{\vartheta}(\rho,\alpha,\beta)} \\
        &\hspace*{3.25cm}d\alpha\,d\beta_1\dots d\beta_{\k-2}
\end{split}\end{equation}
whereas $(\rho,\alpha,\beta)$ denote ($\k$ dimensional) polar coordinates
\begin{gather}
   \vartheta_1=\rho\cos\beta_{\k-2}\dots\cos\beta_1\cos\alpha\notag\\
   \vartheta_2=\rho\cos\beta_{\k-2}\dots\cos\beta_1\sin\alpha\notag\\
   \dots\\
   \vartheta_{\k-1}=\cos\beta_{\k-2}\sin\beta_{\k-3}\notag\\
   \vartheta_{\k}=\sin\beta_{\k-2}\notag
\end{gather}
The factor $\oneover{\k!}$ removes the ordering condition on the simplex,
such that the domain of angle integration is the whole region $[0,\pi/2]^{\k-1}$.
Eventually,
$J_{\vartheta}$ denotes the Jacobi matrix of the coordinate transformation 
$\mapto{(\rho,\alpha,\beta)}{\vartheta}$.
\section{The local equivalence of $\d$ and $\r$ in $\vkn$}
\label{app.drequivalence}
In this appendix the proof of Proposition \ref{prop.localequivdgv} will be
carried out.
Let us recall, what we want to show. Given $\Phi,\Psi$ in the complex
Stiefel manifold $\vkn\subset\CC^{\n\times\k}$, 
the topological distance $\d$ motivated from coding theory is given as
$\d = \fnorm{\Phi-\Psi}$ (we drop the upper index '$V$' in this appendix).

At the same time, locally there is a unique geodesic $\gamma$ in \vkn joining
$\Phi$ and $\Psi$, and the geodesic distance \r is simply defined as its
length $L=\int_0^1 \fnorm{\dot{\gamma}(t)}dt$, $\dot{\gamma}(t)$ being the
parallel transported horizontal tangent vector $X(\gamma(t))$ along
$\gamma$. Thus we obtain $r = \fnorm{X(\gamma(0))}$.
Since both \d and \r are invariant under the action of the \Un we can set
$\Psi=\kbein$ without loss of generality. Recalling the general form 
$X=\textmatrix{A & -B^{\dagger}\\ B & \Zero}$,
$A\in\lieu(\k)$, $B\in\CC^{(\n-\k)\times\k}$, 
of horizontal tangent vectors in $\lieu(\n)$ \eqref{e.Vhorizontal} we arrive at
\begin{gather}
   \d^2 = \fnorm{\Phi-\kbein}^2\\
   \r^2 = \oneover{2}\fnorm{X}^2=\oneover{2}\fnorm{A}^2+\fnorm{B}^2
\end{gather}
whereas $\Phi=\exp X\kbein$.
Unlike the case $A=0$ (representing  tangents for \gkn) there is
no closed form expression for $\Phi$ in terms of $X$ in general 
(compare \cite{ede.ari.smi}), so it remains
a non-trivial task to find constants $\alpha,\beta>0$ satisfying
\begin{equation}
   \boxed{\beta\d\le\r\le\alpha\d} \tag{\ref{e.vchord-vs-geodesic}'}
\end{equation}
expressing the equivalence of \d and \r.

To begin with the easy cases, the constant $\beta$ is easily found, as well
as $\alpha$ when $\k=\n$: Both are simple consequences of the two sided
inequality $\sin x\le x\le (\pi/2)\sin x$, whereas $x\in [0,\pi/2]$ is
understood in the second inequality. 
\begin{lemma}\label{lem.dler}\textin
   In \vkn $\oneover{\sqrt{2}}\d\le\r$ always holds, 
   thus we have $\beta=\oneover{\sqrt{2}}$.
\end{lemma}
\textout
\begin{proof}\myline
   Since $X\in\lieu(\n)$ there exist $V\in\Un$ such that 
   $X=V\diag(\imath\xi)V^{\dagger}$, thus $\r^2=\oneover{2}\norm{\xi}^2$, 
   $\xi=(\xi_1,\dots,\xi_{\n})\in\RR^{\n}$. Now we can estimate as follows
   \begin{equation*}\begin{split}
      \d^2 & = \fnorm{(\Unity-\exp X)\kbein}^2
           \le \fnorm{\Unity-\exp X}^2 \\
           & = \fnorm{\Unity-\exp(\diag(\imath\xi))}^2
             = \sum_j \abs{1-e^{\imath\xi_j}}^2\\
           & = 2\sum_j (1-\cos \xi_j)
             = 4\sum_j \sin^2\frac{\xi_j}{2}
            \le \norm{\xi}^2 = 2\r^2
   \end{split}\end{equation*}
(since $\sin^2 x/2 \le x^2/4$)
\end{proof}
\begin{lemma}\label{lem.dgerkeqn}\textin
   If $\k=\n$ then $\r\le\frac{\pi}{2\sqrt{2}}\d$ holds, 
   thus $\alpha=\frac{\pi}{2\sqrt{2}}$.
\end{lemma}
\textout
\begin{proof}\myline
   $\k=\n$ implies $B=0$, $X=A$ and we can estimate
   \begin{equation*}
      \d^2 = \fnorm{\Unity-\exp A}^2
           = 4\sum_j \sin^2\frac{a_j}{2}
         \ge \frac{4}{\pi^2}\norm{a}^2= \frac{8}{\pi^2}\r^2
   \end{equation*}
   (since $x^2/4\le(\pi^2/4)\sin^2 x/2$ for $x\in [-\pi,\pi]$), 
   whereas $\imath a=\imath(a_1,\dots,a_{\n})$ denotes the vector of eigenvalues of
   $A\in\lieu(\n)$. 
\end{proof}

The non-trivial task is to obtain some $\alpha>0$, when $\k<\n$. The rest
of this section deals with this job.
Let us assume $\k\le\frac{n}{2}$ since this is the relevant case for the
analysis in this work (the case $\k>\frac{n}{2}$ should be similar).
 
Let $X=Y+Z$ with 
$Y=\textmatrix{A & \Zero\\ \Zero & \Zero}$, and
$Z=\textmatrix{\Zero & -B^{\dagger}\\ B & \Zero}$, then 
we can write 
\begin{equation}\label{e.expX-factorised}
   \bar{\Phi}=\exp X=(\exp Z)\textmatrix{v & \Zero\\ \Zero & \Unity}
\end{equation}
since this is merely a factorization of $\Phi=\bar{\Phi}\kbein$ into a
certain projection onto \gkn and the remaining 'phase' in $\Uk\ni v$.
The first factor $\exp Z$ can be calculated in closed form: $B$ has a
singular value decomposition $B=V\diag(\vartheta^{\downarrow})u$
for some $V\in U(\n-\k)$, $u\in\Uk$ and 
$\vartheta^{\downarrow}\mdef(\vartheta_{\k},\dots,\vartheta_1)$ denotes the
vector of principal angles (in decreasing order) between $\generate{\kbein}$ and
$\generate{\Phi}$. 
Setting 
$U_{\vartheta^{\downarrow}}=
 \textmatrix{\diag(\cos\vartheta^{\downarrow}) 
          & -\diag(\sin\vartheta^{\downarrow})
          & \Zero \\
            \diag(\sin\vartheta^{\downarrow})
          & \diag(\cos\vartheta^{\downarrow})
          & \Zero \\
            \Zero & \Zero & \Unity}$,
we arrive at 
$\exp Z = \textmatrix{u & \Zero\\\Zero & V}
          U_{\vartheta^{\downarrow}}
          \textmatrix{u & \Zero\\\Zero & V}^{\dagger}$.
So we have achieved a quite explicit representation of $\bar{\Phi}$. In
particular the principal $\k\times\k$-submatrix 
$\phi=\textmatrix{\Unity & \Zero}\bar{\Phi}\kbein$ reads
\begin{equation}
   \phi = u \diag(\cos\vartheta^{\downarrow})u^{\dagger}v
\end{equation}
Now we can start estimating:
\begin{equation}\begin{split}
   \d^2 & = \fnorm{\Phi-\kbein}^2
          \underset{\fnorm{\Phi}^2=k}{=} 2 (\k-\re\tr\phi) \\
        & = 2 \Big(
                \k-\sum_{j=1}^{\k}
                  \re(u^{\dagger}vu)_{jj}\cos\vartheta_{\k-j}
              \Big) \\
        &\ge 2\Big(
                \k
                -\oneover{2}\sum_j\big[ \re(u^{\dagger}vu)_{jj} \big]^2  
                -\oneover{2}\sum_j \cos^2\vartheta_j
              \Big)
\end{split}\end{equation}
Writing $\Uk\ni v=\exp\tilde{A}$, $\tilde{A}\in\lieu(\k)$ with eigenvalues
$\imath\tilde{a}=\imath(\tilde{a}_1,\dots\tilde{a}_{\k})$ of $\tilde{A}$ 
we have $\re(u^{\dagger}vu)_{jj}\ge 0$ whenever $\tilde{a}\in [-\pi/2,\pi/2]^{\k}$.
Demanding this mild locality restriction we get
$[\re(u^{\dagger}vu)_{jj}]^2 \le \re(u^{\dagger}vu)_{jj}$, thus
$\sum_j [\re(u^{\dagger}vu)_{jj}]^2 
 \le \re\tr(u^{\dagger}vu)=\re\tr v=\sum_j \cos\tilde{a}_j$
and therefore
\begin{equation}\label{e.dtildea}
   \d^2 \ge 2 \sum_j \sin^2\frac{\tilde{a}_j}{2}
           +  \sum_j \sin^2\vartheta_j
        \ge \frac{2}{\pi^2}\norm{\tilde{a}}^2 
           + \frac{4}{\pi^2}\norm{\vartheta}^2
\end{equation}
All what remains to do in order to compare \d with \r is to find the link
between $\tilde{A}$ and $A$, respectively 
$\tilde{Y}=\textmatrix{\tilde{A} & \Zero\\ \Zero & \Zero}$ and 
$Y=\textmatrix{A & \Zero\\ \Zero & \Zero}$. 
By \eqref{e.expX-factorised} 
\begin{equation}
   \exp \tilde{Y} = \eqmatrix{v & \Zero\\ \Zero & \Unity}
                 = \exp(-Z) \exp X
\end{equation}
holds, thus our 'missing link' is given by the Baker-Campbell-Hausdorff
formula expressing $W\in\un$ given by
$\exp W = \exp U \exp V$, $(U,V)\in\un_e\times\un_e$ by
\begin{equation}\label{e.BCH}\begin{split}
   W = & V + \int_0^1 f(e^{t\ad_U} e^{\ad_V}) U dt \\
     = & U + V +
         \sum_{r=1}^{\infty} \frac{(-1)^r}{r+1} \times \\[1ex]
       & \mspace{-18mu}\sum_{\substack{p_1,\dots,p_r\ge 0\\
                           q_1,\dots,q_r\ge 0\\
                           \forall_{i=1..r}\,p_i+q_i>0}}
           \mspace{-18mu}
           \frac{
             \left(
             \frac{\ad_U^{p_1}}{p_1!}\comp\frac{\ad_V^{q_1}}{q_1!}
             \comp\dots\comp
             \frac{\ad_U^{p_r}}{p_r!}\comp\frac{\ad_V^{q_r}}{q_r!}
             \right) (U)
                }
                {p_1+\dots+p_r+1}
\end{split}\end{equation}
whereas $f(z)=\frac{\ln z}{z-1}$ and $\ad_U:\mapto{V}{[U,V]=UV-VU}$
(see \cite{dui.kol} for that particular representation of the BCH formula
(Dynkin's formula in their terminology)). The second part of \eqref{e.BCH} is
nothing but the term-wise integrated Taylor series expansion of the
integrand. Following \cite{dui.kol} the domain of definition $\un_e$ is the
region of \un in which 
the tangent map of $\exp$ is regular. It is the complement of 
$\big\{U\in\un\,|\, 
 \det(\ad_U-2\pi\imath\ZZ'\Unity)=0,\,\ZZ'=\ZZ\setminus\{0\}\big\}$ 
in \un. In particular, $\un_e$ contains a connected neighborhood 
\begin{equation}\label{e.BCH-neighbourhood}
   D(\delta_0)=\{U\in\un\,|\, \fnorm{U}\le\delta_0\}
\end{equation} 
of $\Zero$.
Specializing to $W=\tilde{Y}$, $U=-Z$, $V=X=Y+Z$ yields in multi-index notation
(thus $\abs{p}=\sum_i p_i$, $p!=\prod_i p_i!$)
\begin{equation}\label{e.myBCH}\begin{split}
   \tilde{Y} 
   = & Y + \sum_{r=1}^{\infty} \frac{(-1)^r}{r+1} \times \\
     & \mspace{-24mu}
       \sum_{\substack{p=(p_1,\dots,p_r)\ge 0\\
                       q=(q_1,\dots,q_r)\ge 0\\
                           p+q>0}}
             \mspace{-24mu}          
             \frac{(-1)^{\abs{p}+1}
             \left(
               \ad_Z^{p_1}\comp\ad_X^{q_1}
               \comp\dots\comp
               \ad_Z^{p_r}\comp\ad_X^{q_r}
             \right) (Z)}   
             {(\abs{p}+1)p!q!}
\end{split}\end{equation}
Note that every term contributes at least some factor involving $A$ 
(since $\ad_Z(Z)=0$), hence in the norm estimate 
\begin{equation}\begin{split}
  \fnorm{
    \ad_Z^{p_1} \comp & \ad_X^{q_1}\comp\dots\comp\ad_Z^{p_r}\comp\ad_X^{q_r}(Z)
  } \\
   \le & 2^{\abs{p}+\abs{q}}\fnorm{Z}^{\abs{p}+1}\fnorm{X}^{\abs{q}} \\
   \overset{X=Y+Z}{\le} & 2^{\abs{p}+\abs{q}}\fnorm{Z}^{\abs{p}+1}
       \sum_{i=0}^{\abs{q}}\binom{\abs{q}}{i}\fnorm{Y}^i\fnorm{Z}^{\abs{q}-i}
\end{split}\end{equation}
the term corresponding to $i=0$ has no counterpart in \eqref{e.myBCH}
(resp. it is zero in \eqref{e.myBCH} already),
therefore with $\fnorm{X}\le\delta$ (thus $\fnorm{Y},\fnorm{Z}\le\delta$) 
we can factor out one $\fnorm{Y}$ and estimate
\begin{equation}\begin{split}
   \sum_{i=1}^{\abs{q}} & \binom{\abs{q}}{i}\fnorm{Y}^i\fnorm{Z}^{\abs{q}-i} \\
   & = \Big(\sum_{i=1}^{\abs{q}}
            \binom{\abs{q}}{i}\fnorm{Y}^{i-1}\fnorm{Z}^{\abs{q}-i}
       \Big)\fnorm{Y}\\
   & \le \left(\sum_{i=1}^{\abs{q}}
            \binom{\abs{q}}{i}
       \right)\delta^{\abs{q}-1}\fnorm{Y}
     \le 2^{\abs{q}}\delta^{\abs{q}-1}\fnorm{Y}
\end{split}\end{equation}
and the $\k\times\k$ principal submatrix of \eqref{e.myBCH} of our interest
satisfies
\begin{gather}\label{e.AC}
   \tilde{A} = A + C \\
   \fnorm{C} \le \kappa \fnorm{A} \\
   \label{e.kappa}
   \kappa = \sum_{r=1}^{\infty} \frac{1}{r+1}
              \sum_{\substack{p=(p_1,\dots,p_r)\ge 0\\
                       q=(q_1,\dots,q_r)\ge 0\\
                           p+q>0}}
                   \frac{2^{\abs{p}+2\abs{q}}\delta^{\abs{p}+\abs{q}}}
                        {(\abs{p}+1)p!q!}
\end{gather}
It is possible to rewrite \eqref{e.kappa} such that we can prove the
convergence of the multi-series, that is existence of $\kappa$. given some
sub-multi-indices $p_J$, $q_J$ corresponding to some $J\subset\{1,\dots,r\}$
let us set 
$\lambda_{J}\mdef \frac{(2\delta)^{\abs{p_J}}}{(\abs{p_J}+1)p_J!}$ and
$\mu_{J}\mdef \frac{(4\delta)^{\abs{q_J}}}{q_J!}$, then \eqref{e.kappa}
equals $\sum_r\frac{1}{r+1}\kappa_r$ with 
($J'$ denotes the set $\{1,\dots,r\}\setminus J$)
\begin{equation}\label{e.kappardef}\begin{split}
   \kappa_r 
     = & \sum_{s=0}^{r} 
         \sum_{\substack{J\subset\{1,\dots,r\}\\\abs{J}=s}} \\
       & \Bigg\{
           \big( \sum_{\substack{p_{J}\ge 1\\p_{J'}\equiv 0}}\lambda_{p_{J}} \big)
           \big( \sum_{q\ge 1} \mu_q \big)
           +
           \big( \sum_{p\ge 1} \lambda_p \big)
           \big( \sum_{\substack{q_{J}\ge 1\\q_{J'}\equiv 0}}\mu_{q_{J}} \big)
         \Bigg\}
\end{split}\end{equation}
Now we can perform a rather rough estimate on the sums. We have
$\sum_{q_J\ge 1}\mu_{q_J}=(e^{4\delta}-1)^{\abs{J}}$ and
$\sum_{p_J\ge 1}\lambda_{p_J} \le (e^{2\delta}-1)^{\abs{J}}$, therefore
(note that the sums in the brackets in \eqref{e.kappardef} do not depend on
the particular choice of $J\subset\{1,\dots,r\}$ but only on its
cardinality $\abs{J}=s$)
\begin{equation}\label{e.kapparestimate}\begin{split}
   \kappa_r
   & \le \sum_{s=0}^{r} 
         \binom{r}{s}
         \left\{
           (e^{2\delta}-1)^{s}
           (e^{4\delta}-1)^{r}
           +
           (e^{2\delta}-1)^{r}
           (e^{4\delta}-1)^{s}
         \right\} \\
   & \le (e^{4\delta}-1)^{r} e^{2\delta r}
        +(e^{2\delta}-1)^{r} e^{4\delta r} \\
   & = \big[(e^{4\delta}-1) e^{2\delta}\big]^r 
     + \big[(e^{2\delta}-1) e^{4\delta}\big]^r
\end{split}\end{equation}
It is obvious, that we can choose a $\delta\le\delta_0$ sufficiently small,
such that $\kappa_r\le\oneover{(r+1)^{t}}$ for any given $t>0$, ensuring
the convergence of \eqref{e.kappa}. Setting in particular $t=1$ yields 
$\kappa\le\frac{\pi^2}{6}-1<1$ and we obtain from \eqref{e.AC}
\begin{equation}
   \fnorm{\tilde{A}} \ge (1-\kappa)\fnorm{A}
\end{equation}
Now we can proceed further with \eqref{e.dtildea}:
\begin{equation}\begin{split}
   \d^2 
     & \ge \frac{2(1-\kappa)^2}{\pi^2}\norm{a}^2 
       + \frac{4}{\pi^2}\norm{\vartheta}^2 \\
     & \ge \frac{4(1-\kappa)^2}{\pi^2}\big(\oneover{2}\fnorm{A}^2+\fnorm{B}^2\big)
       = \frac{4(1-\kappa)^2}{\pi^2}\r^2
\end{split}\end{equation}
and we have proven our final lemma:
\begin{lemma}\label{lem.dger}\textin
   If $\k\le\frac{\n}{2}$ there exists a $\delta<\delta_0$, such that
   \mbox{$(1-\kappa)>0$}, whereas $\delta_0$ and $\delta$ are determined by
   \eqref{e.BCH-neighbourhood}, resp. \eqref{e.kapparestimate} demanding 
   $\kappa_r\le\oneover{(r+1)}$.
   Then locally for $r=\fnorm{X}\le\delta$ the relation
   $\r\le\frac{\pi}{2(1-\kappa)}\d$ holds, 
   thus $\alpha=\frac{\pi}{2(1-\kappa)}$.
\end{lemma}
This lemma fills the gap in formula (\ref{e.vchord-vs-geodesic}). Of
course, $(1-\kappa)\approx 1$ would be optimal in this situation (observe
the loss compared to $\alpha$ in Lemma \ref{lem.dgerkeqn}), which can be achieved by
setting $\delta\ll 1$, with $\kappa$ decreasing the smaller
$\delta$ has been chosen. Unfortunately, the smaller we choose 
$\delta$, the higher the required corresponding rate $R$ ensuring the validity of
Lemma \ref{lem.dger}. For example, to obtain a numerical value of 
$R_0\approx 1.4$ 
(by formula (\ref{e.coarsebounds}) as a lower bound for the corresponding
rate, with $\delta=\ubar{\r}_0$), 
which is still achievable for coding purposes in a practical setting, one
needs values of $\delta\approx 1.25$, which is quite large in order to
apply Lemma \ref{lem.dger}, thus the estimates done here are far to rough
to accomplish that. The importance of the lemma actually lies in the fact, that it
proves the \emph{existence} of some $\alpha>0$ in \eqref{e.vchord-vs-geodesic}. 
However, numerical simulations indicate that the real world behaves much
better than the estimates. The histograms in Fig. \ref{fig.hist1-kappa}
display $1-\kappa$ drawn from $1000$ random samples in $\Vc{2,n}$,
$n=4,6,8$ for $\delta=1.25$: 
\begin{figure}[h!tb]
\begin{center}
   \epsfig{file=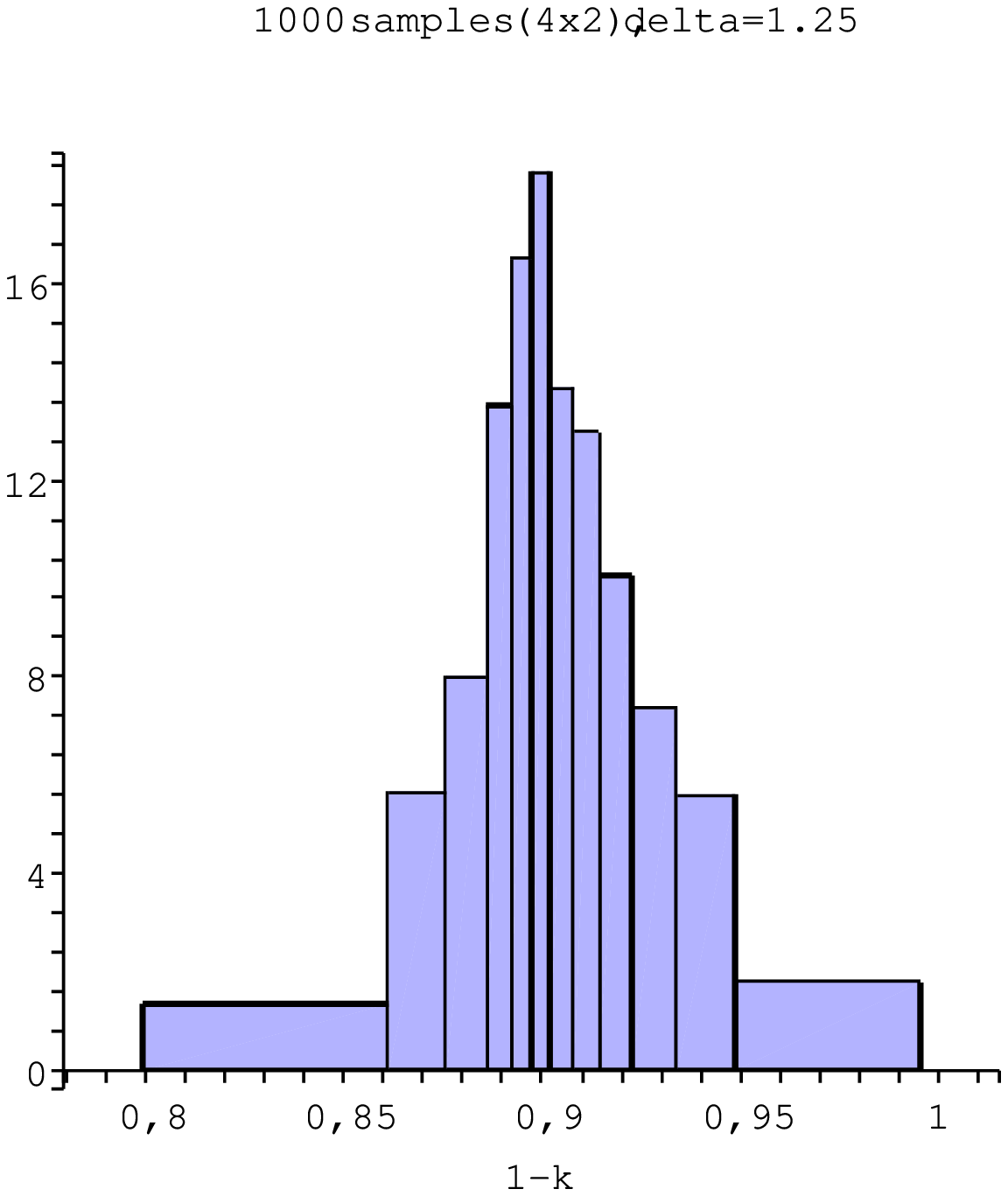,scale=0.45}
   \epsfig{file=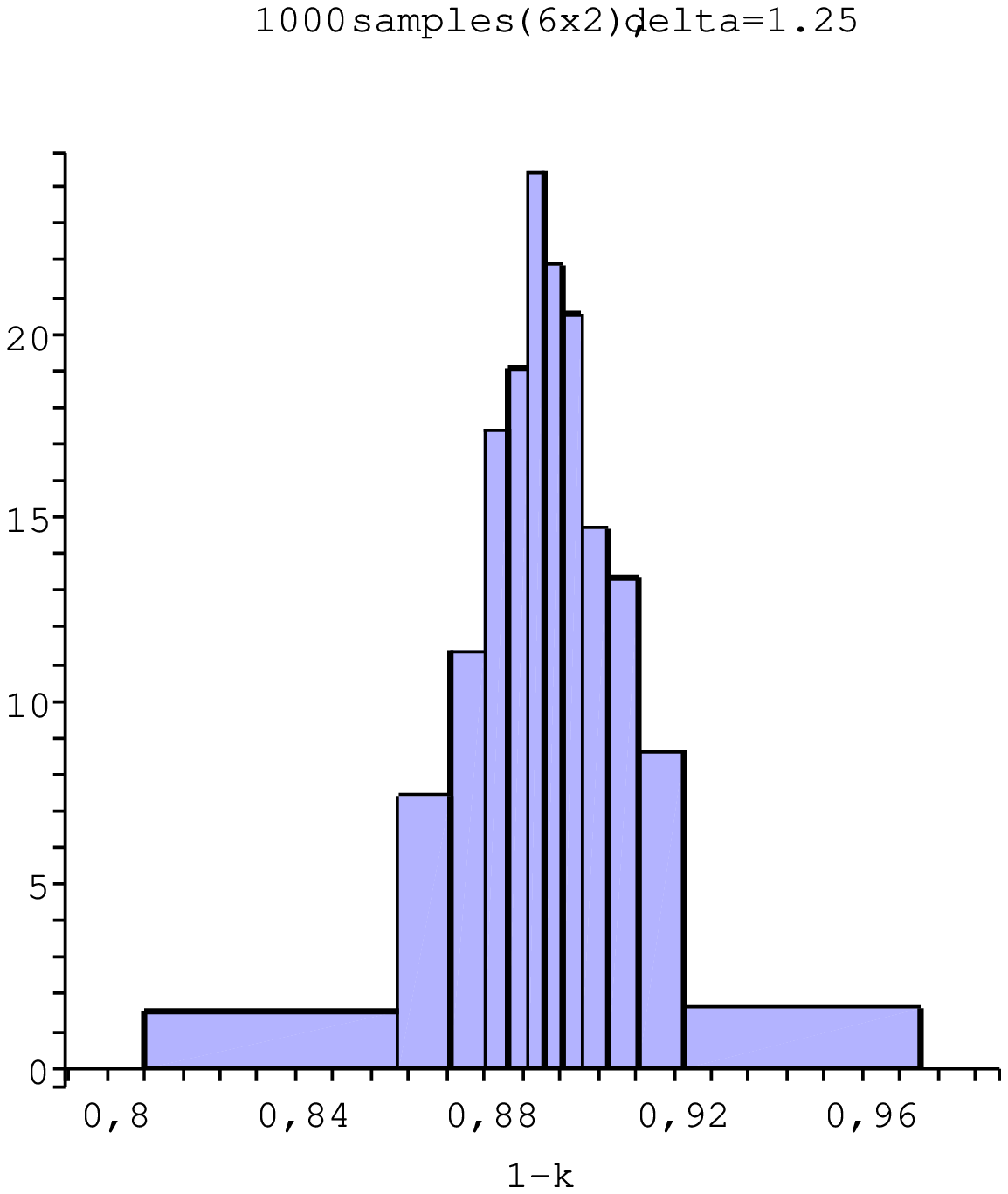,scale=0.45}
   \epsfig{file=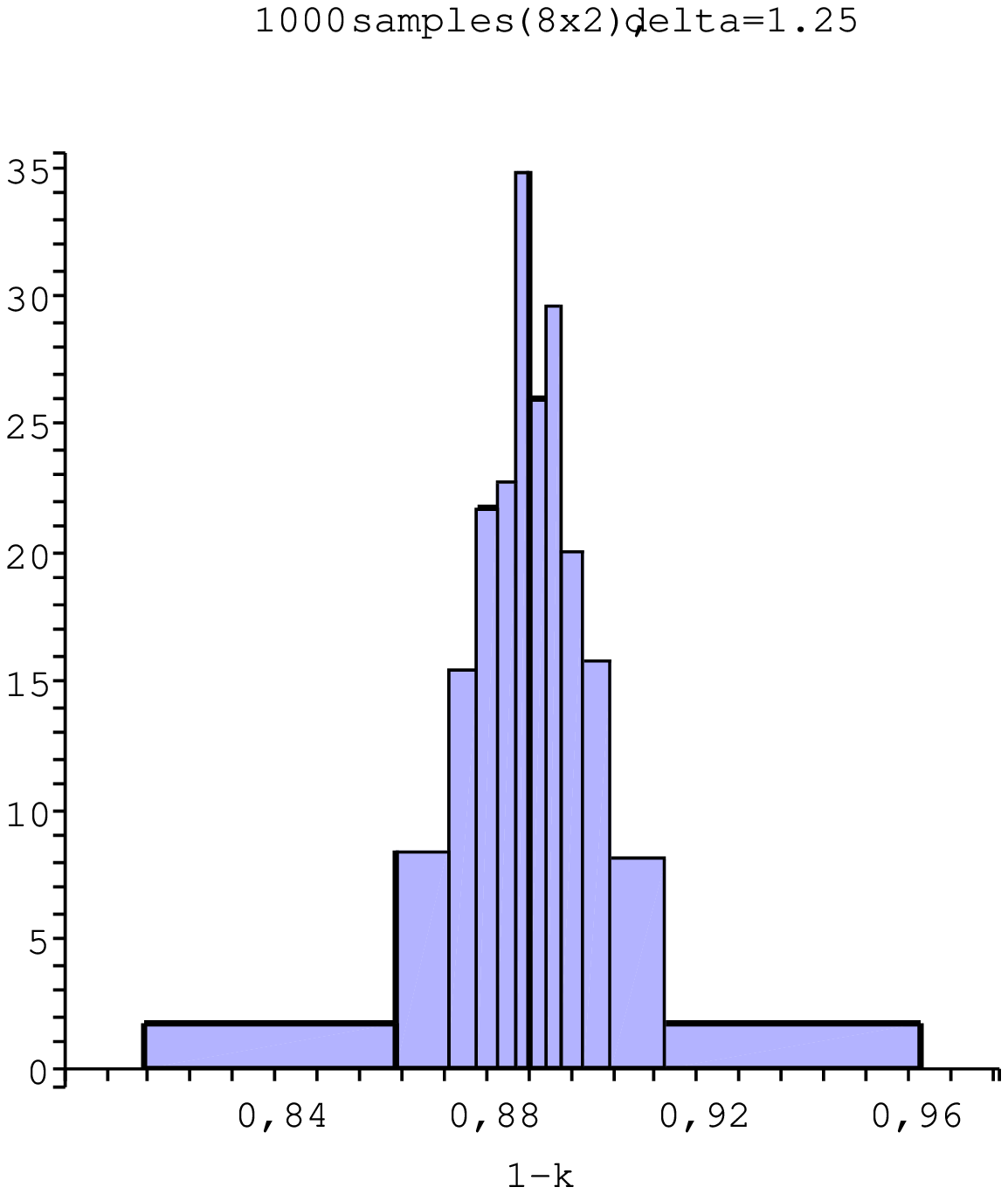,scale=0.45}   
\end{center}\vspace{-2em}
\caption{\label{fig.hist1-kappa}}
\end{figure}
Although there seems to be no rigorous and essentially sharper estimate
available than the one performed here, the numerical examples indicate,
that under still moderate rate constraints $1-\kappa\approx 0.9$ holds, thus 
$\alpha=\frac{\pi}{\sqrt{2}(1-\kappa)}\approx\frac{\pi}{0.9\sqrt{2}}$ in 
\eqref{e.vchord-vs-geodesic}.
%
%
%
%
%
%
%
%
%
%
%
%
\bibliography{../../../../../refhhi}

\begin{thebibliography}{10}
\providecommand{\url}[1]{#1}
\csname url@rmstyle\endcsname
\providecommand{\newblock}{\relax}
\providecommand{\bibinfo}[2]{#2}
\providecommand\BIBentrySTDinterwordspacing{\spaceskip=0pt\relax}
\providecommand\BIBentryALTinterwordstretchfactor{4}
\providecommand\BIBentryALTinterwordspacing{\spaceskip=\fontdimen2\font plus
\BIBentryALTinterwordstretchfactor\fontdimen3\font minus
  \fontdimen4\font\relax}
\providecommand\BIBforeignlanguage[2]{{%
\expandafter\ifx\csname l@#1\endcsname\relax
\typeout{** WARNING: IEEEtran.bst: No hyphenation pattern has been}%
\typeout{** loaded for the language `#1'. Using the pattern for}%
\typeout{** the default language instead.}%
\else
\language=\csname l@#1\endcsname
\fi
#2}}

\bibitem{bar.nog}
A.~Barg and D.~Nogin, ``{Bounds on packings of spheres in the Grassmann
  manifold},'' \emph{{IEEE} Trans. Inform. Theory}, vol.~48, pp. 2450--2454,
  2002.

\bibitem{han.ros}
G.~Han and J.~Rosenthal, ``Unitary space time constellation analysis: An upper
  bound for the diversity,'' 2004, preprint arXiv:math.CO/0401045.

\bibitem{hoc.mar.1}
B.~Hochwald and T.~Marzetta, ``Capacity of a mobile multiple-antenna
  communication link in {Rayleigh} flat fading,'' \emph{{IEEE} Trans. Inform.
  Theory}, vol.~45, pp. 139--157, 1999.

\bibitem{hoc.mar.2}
------, ``Unitary space-time modulation for multiple-antenna communications in
  {Rayleigh} flat fading,'' \emph{{IEEE} Trans. Inform. Theory}, vol.~46, pp.
  543--565, 2000.

\bibitem{zhe.tse}
L.~Zheng and D.~Tse, ``Communication on the {Grassmann} manifold: a geometric
  approach to the noncoherent multiple-antenna channel,'' \emph{{IEEE} Trans.
  Inform. Theory}, vol.~48, pp. 359--383, 2002.

\bibitem{tar.jaf.cal}
V.~Tarokh, H.~Jafarkhani, and A.~Calderbank, ``Space-time block codes from
  orthogonal designs,'' \emph{{IEEE} Trans. Inform. Theory}, vol.~45, pp.
  1456--1467, 1999.

\bibitem{boe.bor.pau}
H.~B{\"o}lcskei, M.~Borgmann, and A.~Paulraj, ``Space-frequency coded
  {MIMO-OFDM} with variable multiplexing-diversity tradeoff,'' \emph{submitted
  to IEEE Trans. Inform. Theory}, 2004.

\bibitem{boe.bor}
H.~B{\"o}lcskei and M.~Borgmann, ``Code design for non-coherent {MIMO-OFDM}
  systems,'' \emph{submitted to IEEE Trans. Inform. Theory}, 2004.

\bibitem{ede.ari.smi}
A.~Edelman, T.~Arias, and S.~Smith, ``{The geometry of algorithms with
  orthogonality constraints},'' \emph{SIAM J. Matrix Anal. Appl.}, vol.~20, pp.
  303--353, 1998.

\bibitem{kue}
W.~K{\"u}hnel, \emph{Differential Geometry: Curves -- Surfaces -- Manifolds},
  ser. Student Mathematical Library, vol 16.\hskip 1em plus 0.5em minus
  0.4em\relax American Mathematical Society, 2002, {German edition:
  Differentialgeometrie: Kurven -- Fl{\"a}chen -- Mannigfaltigkeiten, Vieweg
  1999}.

\bibitem{con.har.slo}
J.~Conway, R.~Hardin, and N.~Sloane, ``{Packing lines, planes, etc.: Packings
  in Grassmannian spaces},'' \emph{Experimental Mathematics}, vol.~5, pp.
  139--159, 1996, url: http://www.research.att.com/{\raisebox{-0.75ex}{\~{
  }}}njas/grass/index.html.

\bibitem{gal.hul.laf}
S.~Gallot, D.~Hulin, and J.~Lafontaine, \emph{Riemannian Geometry},
  2nd~ed.\hskip 1em plus 0.5em minus 0.4em\relax Springer, 1993.

\bibitem{lia.xia}
X.-B. Liang and X.-G. Xia, ``Unitary signal constellations for differential
  space-time modulation with two transmit antennas: Parametric codes, optimal
  designs, and bounds,'' \emph{{IEEE} Trans. Inform. Theory}, vol.~48, pp.
  2291--2322, 2002.

\bibitem{foz.mclau.schaf}
M.~Fozunbal, S.~McLaughlin, and R.~Schafer, ``On performance limits of
  {Space-time} codes: A sphere-packing bound approach,'' \emph{{IEEE} Trans.
  Inform. Theory}, vol.~49, no.~10, 2003.

\bibitem{poor}
W.~A. Poor, \emph{Differential Geometric Structures}.\hskip 1em plus 0.5em
  minus 0.4em\relax McGraw-Hill Inc., 1981.

\bibitem{has.hoc.mar}
B.~Hassibi, B.~Hochwald, and T.~Marzetta, ``Space-time autocoding,''
  \emph{{IEEE} Trans. Inform. Theory}, vol.~47, pp. 2761--2781, Nov 2001.

\bibitem{dui.kol}
J.~Duistermaat and J.~Kolk, \emph{Lie Groups}.\hskip 1em plus 0.5em minus
  0.4em\relax Springer, 2000.

\end{thebibliography}
%
%
%
%
%
%
\end{document}